\documentclass[autodetect-engine,dvipdfmx-if-dvi,ja=standard, 12pt, a4paper]{article}
\usepackage{amsmath}
\usepackage{ascmac}
\usepackage{amssymb}
\usepackage{graphicx}
\usepackage{longtable}
\usepackage{bm}
\usepackage{enumitem}
\newcommand{\ve}{\bm{e}} 

\newcommand{\vv}{{\bm{v}}} 
\newcommand{\vw}{{\bm{w}}} 

\newcommand{\vH}{\bm{H}}
\newcommand{\vX}{\bm{X}}

\begin{document}

\title{Mean curvature flow for principal orbits of\\
Hermann actions on rank two symmetric spaces}
\author{Naoyuki Koike and Sakura Nakaoka}
\date{}
\maketitle

\begin{abstract}
In this paper, we illustrate the behaviour of the mean curvature flows starting from principal orbits of 
any commuting Hermann action of cohomogeneity two on irreducible rank two Riemannian symmetric spaces of compact type 
by using Mathematicae.  In more detail, we illustrate the velocity vector fields of the curves (determined by the flows) 
on the orbit space (which is a 2-simplex) of the Hermann action by using Mathematcae.  
Also, we calculate the position of the point of the orbit space corresponding to the only minimal principal orbit 
of the Hermann action by using Mathematicae.  
\end{abstract}





\section*{$\S0.$ Introduction}
\addcontentsline{toc}{part}{$\S0.$ Introduction}
An isometric action of a compact Lie group $G$ on a Riemannian manifold $M$ is said to be {\it polar} if there exists 
a closed, connected submanifold $\Sigma$ of $M$ meeting all $G$-orbits orthogonally.  The submanifold $\Sigma$ is called 
a {\it section} of the action.  A section is a automatically totally geodesic in $M$.  Furthermore, if the induced metric 
on the section is flat, then the action is said to be {\it hyperpolar}.  
The notion of a hyperpolar action was introduced by Heintze-Palais-Terng-Thorbergsson (\cite{HPTT}).  
If $(G,K)$ and $(G,H)$ are symmetric pairs 
(i.e, $G/K$ and $G/H$ are Riemanian symmetric spaces of compact type), then the natural action $H\curvearrowright G/K$ 
is called a {\it Hermann action}.  Let $\theta$ and $\tau$ of $G$ be an involutions of $G$ satisfying 
$({\rm Fix}\,\theta)_0\subset K\subset {\rm Fix}\,\theta$ and $({\rm Fix}\,\tau)_0\subset H\subset{\rm Fix}\,\tau$, 
respectively, where ${\rm Fix}\,(\bullet)$ is the fixed point group of $\bullet$ and $({\rm Fix}(\bullet))_0$ 
is the identity component of ${\rm Fix}(\bullet)$.  In particular, if $\theta$ and $\tau$ commute to each other, 
then this action is called a {\it commuting Hermann action}.    
It is known that Hermann actions are hyperpolar ([HPTT]).  Also, it is known that all hyperpolar actions of cohomogeneity 
greater than one on irreducible Riemannian symmetric spaces of compact type are Hermann ones ([Kol]).  

As a generalized notion of compact isoparametric hypersurfaces in a sphere and a hyperbolic space, 
and compact isoparametric submanifolds in a Euclidean space, Terng-Thorbergsson [TT] introduced the notion of 
an equifocal submanifold in a Riemannian symmetric space as a compact submanifold M satisfying the following 
three conditions:\\

\vspace{0.15truecm}

\noindent
(i)\ The normal holonomy group of $M$ is trivial;

\noindent
(ii)\ $M$ has a flat section, that is, for each $x \in M$, $\Sigma_x :=$ exp$^{\perp}(T_x^{\perp}M)$ is totally geodesic and\\
\hspace{15pt} the induced metric on $\Sigma_x$ is flat, where $T_x^{\perp}M$ is the normal space of $M$ at $x$ and exp$^{\perp}$ is 
\\\hspace{15pt} the normal exponential map of $M$;

\noindent
(iii)\ For each parallel normal vector field $v$ pf $M$, the focal radii of $M$ along the normal geode-\\\hspace{15pt} 
sic $\gamma_{v_x}$ (with $\gamma'_{v_x}(0) = v_x$) are independent of the choice of $x \in M$ with considering \\\hspace{15pt}
the multiplicities, where $\gamma'_{v_x}(0)$ is the velocity\\\hspace{15pt} vector of $\gamma_{v_x}$ at 0.  

\vspace{0.15truecm}

\noindent
According to results in [C], [HPTT] and [Kol], all equifocal submanifolds of cohomogeneity greater than one 
in irreducible Riemannian symmetric spaces of compact type occur as principal orbits of Hermann actions 
on the Riemannian symmetric spaces.\\

Let $\{f_t\}_{t\in [0,T)}$ be a $C^{\infty}$-family of $C^{\infty}$-immersions of a manifold $M$ into 
a Riemannian manifold $N$, where $T$ is a positive constant or $T=\infty$. Define a map 
$\widetilde{f}:M \times [0,T) \rightarrow N$ by $\widetilde{f}(x,t) = f_t(x)\ ((x,t) \in M \times [0,T))$. 
If, for each $t \in [0,T)$ and any $x\in M$, $\widetilde{f}_*((\frac{\partial}{\partial t})_{(x,t)})$ is equal to 
the mean curvature vector of $f_t : M \hookrightarrow N$ at $x$, then $\{f_t\}_{t \in [0,T)}$ is called 
a {\it mean curvature flow}.  In particular, if $f_t$'s are embeddings, then we call $\{M_t := f_t(M)\}_{t\in [0,T)}$ 
rather than $\{f_t\}_{t\in [0,T)}$ a mean curvature flow.  

The mean curvature flows starting from isoparametric submanifolds and their focal submanifolds in spheres were studied 
in \cite{LT}.  Later, the mean curvature flows for equifocal submanifolds and their focal submanifolds 
in Riemannian symmetric spaces of compact type were studied in \cite{Koi}.  In particular, the mean curvature flows 
starting from the principal orbits and singular ones of commuting Hermann actions of cohomogeneity two 
on irreducible rank two Riemannian symmetric spaces of compact type were studied in detail in \cite{Koi}.  
The main result in \cite{Koi} is as follows.\\
\\
\textbf{FACT A} (\cite{Koi})\ \ {\sl
Let $M$ be an equifocal submanifold in a Riemannian symmetric space $G/K$ of compact type. Then the following statements {\rm (i)} and {\rm (ii)} hold :
\begin{description}[labelwidth=24pt]
\item[{\rm (i)}] If $M$ is not minimal, then the mean curvature flow $M_t$ having $M$ as initial data collapses to a focal submanifold $F$ of $M$ in finite time. Furtheremore, if $M$ is irreducible, the codimention of $M$ is greater than one and if the fibration of $M$ onto $F$ is spherical, then the flow $M_t$ has type I singularity.
\item[{\rm (ii)}] For any focal submanifold $F$ of $M$, there exists a parallel submanifold of $M$ collapsing to $F$ along the mean curvature flow and the set of all parallel submanifolds collapsing to $F$ along the mean curvature flow is a one-parameter $C^{\infty}$-family.
\end{description}}
\vspace{20pt}

\section*{$\S1.$ Basic notions and facts}
\addcontentsline{toc}{part}{$\S1.$ Basic notions and facts}

In this section, we recall the notions of a Riemannian symmetric space of compact type and commuting Hermann actions, and 
basic facts related to these notions.  

\subsection*{$\S\S1.1.$ The restricted root system of a symmetric pair}
\addcontentsline{toc}{section}{$\S\S1.1.$ The restricted root system of a symmetric pair}

Let $G$ be a compact semi-simple connected Lie group, and $\mathfrak{g}$ its Lie algebra, which is identified with 
the tangent space $G$ at the identity element of $G$.  A subgroup $K$ of $G$ is called a {\it symmetric subgroup} of $G$ 
(and the pair $(G,K)$ is called a {\it symmetric pair}) if there is an involution $\theta$ of $G$ such that 
$({\rm Fix}\,\theta)_0 \subset G \subset{\rm Fix}\,\theta$.  In this case, the quotient manifold $G/K$ equipped with 
a $G$-invariant metric induced from any ${\rm Ad}(G)$-invariant inner product of $\mathfrak{g}$ is called a 
{\it Riemannian symmetric space of compact type}.\\

First, we recall the notion of a root system.  Let $(V,\langle\,, \rangle)$ be a $r$-dimensional inner product vector space and $\Delta$ be a subset of $V^*\backslash\{{\bm 0}\}$. Define a reflection $s_{\alpha} : V \rightarrow V \ (\alpha \in \Delta)$ by
\begin{equation*}
s_{\alpha}({\bm v}) := {\bm v} - \frac{2\langle {\bm v}, {\bm v}_{\alpha}\rangle}{\langle {\bm v}_{\alpha}, {\bm v}_{\alpha}\rangle} \cdot {\bm v}_{\alpha} \quad ({\bm v} \in V^*),
\end{equation*}
where ${\bm v}_{\alpha}$ is the vector of $V$ with $\alpha(\cdot) := \langle {\bm v}_{\alpha}, \cdot\rangle$.\\
\\
\textbf{Definition.}\ \ $\Delta$ is called a {\it root system of rank} $r$ if it satisfies the following conditions :
\begin{description}[labelwidth=24pt] 
\item[{\rm (i)}] Span$\{ {\bm v}_{\alpha} \mid \alpha \in \Delta \} = V$ \vspace{-5pt}
\item[{\rm (ii)}] $\dfrac{2\langle {\bm v}_{\alpha}, {\bm v}_{\beta}\rangle}{\langle {\bm v}_{\alpha}, {\bm v}_{\alpha}\rangle}\;(\alpha, \beta \in \Delta)$ are integer numbers.
\end{description}
\vspace{15pt}

\noindent
\textbf{Definition.}\ \ A subsystem $\Pi = \{ \alpha_1, ... , \alpha_r\}$ of $\Delta$ is called a {\it fundamental root system} if it satisfies the following conditions :
\begin{description}[labelwidth=24pt] 
\item[{\rm (i)}] $(\alpha_1, ... , \alpha_r)$ is a base of $V^*$ \vspace{-5pt}
\item[{\rm (ii)}] Each element $\alpha$ of $\Delta$ is described as $\alpha = \sum_{i=1}^{r}n_i\alpha_i \; (n_i \geq 0 \;(i = 1, ... , r)$ or $n_i \leq 0 \;(i = 1, ... , r) )$.
\end{description}
\vspace{15pt}

Let $\widetilde{\mathfrak{g}}$ be a complex semi-simple Lie algebra and $\widetilde B$ its killing form.  
This form $\widetilde{B}$ is a non-degenerate complex bilinear form on $\widetilde{\mathfrak{g}}$ defined by 
$\widetilde{B}(\bm{v},\bm{w}) := {\rm Tr}({\rm ad}(\bm{v})\circ{\rm ad}(\bm{w}))$ 
($\vv,\,\vw\in\widetilde{\mathfrak{g}}$).  
A maximal abelian subalgebra $\widetilde{\mathfrak{a}}$ of $\widetilde{\mathfrak{g}}$ is called the {\it Cartan subalgebra}, if a complex linear transformation ad$(\bm{v}) : \widetilde{\mathfrak{g}}\rightarrow\widetilde{\mathfrak{g}}\ (\forall\bm{v}\in \widetilde{\mathfrak{a}})$ is diagonalizable.  
Since the family \{ad$({\bm v})\mid{\bm v}\in\widetilde{\mathfrak{a}}$\} is a commmuting family of skew-symmetric 
complex linear operators, 
they are diagonalizable simultaneously, that is, there is a simultaneously eigenspace decomposition 
$\widetilde{\mathfrak{g}} = \oplus^k_{i=1}\widetilde{\mathfrak{g}}_i$. $\widetilde{\mathfrak{g}}_i$ is given by
\begin{center}
$\widetilde{\mathfrak{g}}_i = \underset{{\bm v} \in \widetilde{\mathfrak{a}}}{\cap}$ Ker(ad$({\bm v}) - \alpha^{{\bm v}}_i$ id$_{\widetilde{\mathfrak{g}}})$
\end{center}
($\alpha^{{\bm v}}_i$ is one of complex eigenvalues of ad$({\bm v})$).  We can show that 
the correpondence $\alpha_i : {\bm v} \mapsto  \alpha^{{\bm v}}_i\ ({\bm v} \in \widetilde{\mathfrak{a}})$ 
is a complex linear function.  
Let $\widetilde{\mathfrak{a}}^*$ be the dual vector space of complex vector space $\widetilde{\mathfrak{a}}$. Define a complex vector subspace $\widetilde{\mathfrak{g}}_{\alpha}$ of $\widetilde{\mathfrak{g}}$ by
\begin{center}
$\widetilde{\mathfrak{g}}_{\alpha} := \{\bm{w}\in\widetilde{\mathfrak{g}}\mid$ ad$(\bm{v})(\bm{w})=\alpha(\bm{v})\cdot\bm{w}\ (\forall\bm{v}\in\widetilde{\mathfrak{a}})\}\ (\alpha\in\widetilde{\mathfrak{a}}^*)$
\end{center}
and a subset $\widetilde{\Delta}$ of $\widetilde{\mathfrak{a}}^*$ by
\begin{center}
$\widetilde{\Delta} := \{\alpha\in\widetilde{\mathfrak{a}}^*\backslash\{\bm{0}\}\mid\widetilde{\mathfrak{g}}_{\alpha}\neq\{\bm{0}\}\}$.
\end{center}
From $\widetilde{\mathfrak{g}} = \oplus^k_{i=1}\widetilde{\mathfrak{g}}_i$, we have the following simultaneously eigenspace decomposition :
\begin{equation}
\widetilde{\mathfrak{g}} = \widetilde{\mathfrak{a}} \oplus \left( \underset{\alpha \in \widetilde{\Delta}}{\oplus}\ \widetilde{\mathfrak{g}}_{\alpha} \right). \tag{1.1.1}
\end{equation}
Define ${\bm v}_{\alpha} \in \widetilde{\mathfrak{a}}\ (\alpha \in \widetilde{\Delta})$ by $\widetilde{B}({\bm v}_{\alpha}, \cdot) = \alpha (\cdot)$.
Let $\widetilde{\mathfrak{a}}_{\mathbb{R}}$ be a real form  of $\widetilde{\mathfrak{a}}$ by
\begin{center}
$\widetilde{\mathfrak{a}}_{\mathbb{R}} :=$ Span$_{\mathbb{R}} \{ {\bm v} \mid \alpha \in \widetilde{\Delta} \}$,
\end{center}
where Span$_{\mathbb{R}}\{\bullet\}$ is the real vector subspace spanned by $\bullet$.
We have $\alpha (\widetilde{\mathfrak{a}}_{\mathbb{R}}) \subset \mathbb{R}\ (\alpha \in \widetilde{\Delta})$ and 
\begin{center}
$\widetilde{\Delta}_{\mathbb{R}} := \{\alpha|_{\widetilde{\mathfrak{a}}_{\mathbb{R}}}\mid \alpha \in \widetilde{\Delta} \}\ (\subset \widetilde{\mathfrak{a}}^*_{\mathbb{R}})$
\end{center}
is a root system where $\widetilde{\mathfrak{a}}^*_{\mathbb{R}}$ is the dual space of the real vector space. 
$\widetilde{\Delta}_{\mathbb{R}}$ is called the {\it root system} of $\widetilde{\mathfrak{g}}$ with respect to $\widetilde{\mathfrak{a}}$, the decomposition (1.1.1) is called the {\it root space decomposition of} $\widetilde{\mathfrak{g}}$ {\it with respect to} 
$\widetilde{\mathfrak{a}}$, and $\widetilde{\mathfrak{g}}_{\alpha}\ (\alpha \in \widetilde{\Delta})$ is called the {\it root space} 
for $\alpha$.

Let $(\mathfrak{g}, \theta)$ be a semi-simple orthogonal symmetric Lie algebra of compact type, $B$ be the killing form of $\mathfrak{g}$. Set $\mathfrak{k}:=$ Ker$(\theta\ -$ id$_{\mathfrak{g}})$, $\mathfrak{p}:=$ Ker$(\theta\ +$ id$_{\mathfrak{g}})$.  
Let $\mathfrak{a}$ be a maximal abelian subspace of $\mathfrak{p}$ and $\widetilde{\mathfrak{a}}$ be a Cartan subalgebla ofthe complexification $\mathfrak{g}^{\mathbb{C}}$ of $\mathfrak g$ including $\mathfrak{a}$.  
The dimension of $\mathfrak{a}$ depends only on $(\mathfrak{g}, \theta)$.  This dimension is called the {\it rank} of 
$(\mathfrak{g}, \theta)$.  Let $\widetilde{\Delta}$ be a root system of $\mathfrak{g}^{\mathbb{C}}$ with respect to $\widetilde{\alpha}$ and $\mathfrak{g}_{\alpha}^{\mathbb{C}}$ be a root space for $\alpha \in \widetilde{\Delta}$. 
Since $(\mathfrak{g}, \theta)$ is of compact type, $\alpha|_{\mathfrak{a}}$ is a purely imaginary-valued function.\\
Define a subset $\Delta$ of $\mathfrak{a}^*$ by
\begin{align*}
\Delta &:= \{ \sqrt{-1} \alpha|_{\mathfrak{a}}\mid\alpha \in \widetilde{\Delta}\ {\rm s.t.}\ \alpha|_{\mathfrak{a}}\neq 0\}.
\end{align*}
For a linear function $\lambda$ on $\mathfrak{a}$, define a vector subspace $\mathfrak{p}_{\lambda}$ of $\mathfrak{p}$ by
\begin{center}
$\mathfrak{p}_{\lambda} := \{\boldsymbol{w} \in \mathfrak{p}\mid$ ad$(\boldsymbol{v})^2(\boldsymbol{w}) = -\lambda(\boldsymbol{v})^2\cdot\boldsymbol{w}\ (\forall\boldsymbol{v}\in\mathfrak{a})\}$
\end{center}
and a vector subspace $\mathfrak{k}_{\lambda}$ of $\mathfrak{k}$ by
\begin{center}
$\mathfrak{k}_{\lambda} := \{\boldsymbol{w} \in \mathfrak{k}\mid$ ad$(\boldsymbol{v})^2(\boldsymbol{w}) = -\lambda(\boldsymbol{v})^2\cdot\boldsymbol{w}\ (\forall\boldsymbol{v}\in\mathfrak{a})\}$.
\end{center}
Denote by $\geq$ the lexicographic order of $\mathfrak{a}^*$ associated to a base of $\mathfrak{a}$. 
Set $\Delta_+ := \{ \lambda \in \Delta \mid \lambda \geq {\bm 0} \}$. 
We call $\Delta$ the {\it restricted root system of} $(\mathfrak{g}, \theta)$ {\it with respect to} $\mathfrak{a}$ 
and $\Delta_+$ the {\it positive restricted root system of} $(\mathfrak{g}, \theta)$ {\it with respect to} 
$\mathfrak{a}$.  Then the following proposition holds.\\
\\
\textbf{Proposition 1.1.3.}\ \ {\sl The following statements (i)-(iii) hold.  
\begin{description}[labelwidth=24pt] 

\item[{\rm (i)}] The system $\Delta$ is a root system.

\item[{\rm (ii)}] For any $\lambda\in\Delta$,
\begin{center}
$\mathfrak{p}_{\lambda} = \left( \underset{\alpha\in\widetilde{\Delta}_{\lambda}}{\oplus} (\mathfrak{g}_{\alpha}^{\mathbb{C}}\oplus\mathfrak{g}_{-\alpha}^{\mathbb{C}})\right)\cap\mathfrak{p}$,  $\mathfrak{k}_{\lambda} = \left( \underset{\alpha\in\widetilde{\Delta}_{\lambda}}{\oplus} (\mathfrak{g}_{\alpha}^{\mathbb{C}}\oplus\mathfrak{g}_{-\alpha}^{\mathbb{C}})\right)\cap\mathfrak{k}$
\end{center}
($\widetilde{\Delta}_{\lambda} := \{\alpha\in\widetilde{\Delta}\mid\alpha|_{\mathfrak{a}} = \pm\sqrt{-1}\lambda\}$).  
Also $\dim\mathfrak{p}_{\lambda} = \dim\mathfrak{k}_{\lambda}$ = $\frac{1}{2}\cdot\sharp\widetilde{\Delta}_{\lambda}$ holds, where $\sharp\widetilde{\Delta}_{\lambda}$ is the cardinal number of $\widetilde{\Delta}_{\lambda}$.

\item[{\rm (iii)}] We have the following simultaneously eigenspace decompositions :
\begin{center}
$\mathfrak{p}$ = $\mathfrak{a} \oplus \left( \underset{\lambda \in \Delta_+}{\oplus} \mathfrak{p}_{\lambda} \right)$\ \ and\ \ 
$\mathfrak{k}$ = $\mathfrak{z}_{\mathfrak{k}}(\mathfrak{a}) \oplus \left( \underset{\lambda \in \Delta_+}{\oplus} \mathfrak{k}_{\lambda} \right),$
\end{center}
where $\mathfrak{z}_{\mathfrak{k}}(\mathfrak{a})$ is the centralizer of $\mathfrak{a}$ in $\mathfrak{k}$.  
\end{description}}
\vspace{20pt}

\subsection*{$\S\S1.2.$ The orbit geometry of Hermann actions}
\addcontentsline{toc}{section}{$\S\S1.2.$ The orbit geometry of Hermann actions}

Let $G/K$ be a Riemannian symmetric space of compact type and $H$ a symmetric subgroup of $G$.  Also, let $\theta$ 
(resp. $\tau$) be an involution of $G$ with $({\rm Fix}\,\theta)_0 \subset \mathit{K} \subset {\rm Fix}\,\theta$ (resp. 
$({\rm Fix}\,\tau)_0 \subset \mathit{H} \subset {\rm Fix}\,\tau$).  
We give $G/K$ the $G$-invariant metric $\widetilde g$ with $\widetilde{g}_{eK}=-B|_{\mathfrak{p}\times\mathfrak{p}}$, 
where $B$ is the killing form of $\mathfrak g$.  
In the sequel, we assume that $[\theta, \tau] = 0$.
Then the natural action $H \curvearrowright G/K$ is called {\it commuting Hermann action}.  
Set $\mathfrak{g}$ := Lie$\,G$, $\mathfrak{k}$ := Lie$\,\mathit{K}$, $\mathfrak{h}$ := Lie$\,\mathit{H}$, $\mathfrak{p}$ := Ker($\theta$ + id), $\mathfrak{q}$ := Ker($\tau$ + id). 
From $\theta \circ \tau = \tau \circ \theta$, we have
\begin{center}
$\mathfrak{p} = (\mathfrak{p} \cap \mathfrak{h}) \oplus (\mathfrak{p} \cap \mathfrak{q}$), $\mathfrak{k} = (\mathfrak{k} \cap \mathfrak{h}) \oplus (\mathfrak{k} \cap \mathfrak{q}$).
\end{center}
Take a maximal abelian subspace $\mathfrak{b}$ of $\mathfrak{p} \cap \mathfrak{q}$ and a maximal abelian subspace $\mathfrak{a}$ of $\mathfrak{p}$ with $\mathfrak{b} \subset \mathfrak{a}$. 
Let $\Delta$ be a restricted root system with respect to $\mathfrak{a}$, and the root space deconposition of $\mathfrak{k}$ and $\mathfrak{p}$ with respect to $\mathfrak{a}$ be
\begin{center}
$\mathfrak{k}$ = $\mathfrak{z}_{\mathfrak{p}}(\mathfrak{a}) \oplus \left( \underset{\lambda \in \Delta_+}{\oplus} \mathfrak{k}_{\lambda} \right)$, $\mathfrak{p}$ = $\mathfrak{a} \oplus \left( \underset{\lambda \in \Delta_+}{\oplus} \mathfrak{p}_{\lambda} \right).$
\end{center}
Define $\Delta'_{V}$ and $\Delta'_{H}$ by
\begin{center}
$\Delta'_{V} := \{\lambda|_{\mathfrak{b}}\mid\lambda\in\Delta$ s.t. $\mathfrak{p}_{\lambda} \cap \mathfrak{q}\neq\{0\}\}$, \\
$\Delta'_{H} := \{\lambda|_{\mathfrak{b}}\mid\lambda\in\Delta$ s.t. $\mathfrak{p}_{\lambda} \cap \mathfrak{h}\neq\{0\}\}$.
\end{center}
Set $(\Delta'_V)_+ := \Delta'_V \cap \Delta_+$ and $(\Delta'_H)_+ := \Delta'_H \cap \Delta_+$.  
Note that $\Delta'_{V}$ is a root system but $\Delta'_{H}$ is not necessarily a root system.  
Since $\mathfrak{p} \cap \mathfrak{q}$ and $\mathfrak{p} \cap \mathfrak{h}$ are ad$_{\mathfrak{g}}({\bm v})^2$-invariant 
for any $\vv\in\mathfrak b$, we have 
\begin{align*}
 \mathfrak{p} \cap \mathfrak{q} &= \mathfrak{b} \oplus \left( \underset{\lambda \in (\Delta'_{V})_+}{\oplus}(\mathfrak{p}_{\lambda} \cap \mathfrak{q})\right),\\
\mathfrak{p} \cap \mathfrak{h} &= \mathfrak{z}_{\mathfrak{p} \cap \mathfrak{h}}(\mathfrak{b}) \oplus \left( \underset{\lambda \in (\Delta'_{H})_+}{\oplus}(\mathfrak{p}_{\lambda} \cap \mathfrak{h})\right).
\end{align*}
Set $L$ := Fix$(\tau \circ \theta)$ and denote by $\mathfrak{l}$ its Lie algebra.  
The tangent space of the orbit $F := H(eK)$ at $eK$ is identified with $\mathfrak{p} \cap \mathfrak{h}$.  
Since $\mathfrak{p} \cap \mathfrak{h}$ is a Lie triple system, $F = {\rm exp}_{eK}(\mathfrak{p} \cap \mathfrak{h})$ is 
totally geodesic in $G/K$ and it is the Riemannian symmetric space of compact type associated to the semi-simple 
symmetric pair $(H, H \cap K)$.  
On the other hand, a tangent space of the normal umbrella $F^{\perp} := {\rm exp}^{\perp}(T^{\perp}_{eK}F)$ of $F$ 
at $eK$ is identified with $\mathfrak{p} \cap \mathfrak{q}$.  
Since $\mathfrak{p} \cap \mathfrak{q}$ is a Lie triple system, $F^\perp = {\rm exp}_{eK}(\mathfrak{p} \cap \mathfrak{q})$ 
is totally geodesic in $G/K$ and it $F^{\perp}$ is the Riemannian symmetric space of compact type associated to 
the semi-simple symmetric pair $(L, H \cap K)$. \\

Define a subset $C'$ of $\mathfrak{b}$ by
\begin{center}
$C':=\left\{ Z \in \mathfrak{b} \mid 0 < \lambda(Z) < \pi\ (\forall \lambda \in (\Delta'_V)_+),\ -\dfrac{\pi}{2} < \lambda(Z) < \dfrac{\pi}{2}\ (\forall \lambda \in (\Delta'_H)_+) \right\}.$
\end{center}
All $H$-orbits intersect with $\exp_{eK}(\overline{C'})$ only at a one point, that is, an orbit space $H \backslash G/K$ 
of $H \curvearrowright G/K$ is identified with $\overline{C'}$, where $\overline{C'}$ denotes the closure of $C'$.  
The following facts hold.\\
\\
\textbf{Proposition 1.2.1.}\ \ {\sl Let $p:= \exp_{eK}(Z) \in \exp_{eK}(\overline{C'})\ (Z \in \overline{C'})$.  
\begin{description}[labelwidth=24pt] 
\item[{\rm (i)}]\ The tangent space $T_p(H\cdot p)$ is decomposed orthogonally as follows:
\begin{align*}
T_p(H \cdot p) 
&= (T_p(H \cdot p) \cap P_{\gamma_{Z}|_{[0,1]}}(T_{eK}F)) \oplus (T_p(H \cdot p) \cap T_p(F^{\perp}))\\
&= P_{\gamma_{Z}|_{[0,1]}} \left(\mathfrak{z}_{\mathfrak{p} \cap \mathfrak{h}}(\mathfrak{b}) \oplus 
\left( \underset{\lambda \in (\Delta'_{H})_{Z}+}{\oplus}(\mathfrak{p}_{\lambda} \cap \mathfrak{h})\right)
\oplus
\left( \underset{\lambda \in (\Delta'_{V})_{Z}}{\oplus}(\mathfrak{p}_{\lambda} \cap \mathfrak{q})\right)\right),
\end{align*}
where $\gamma_{Z}$ denotes the geodesic in $G/K$ with $\gamma_{Z}'(0)=Z$, 
$P_{\gamma_{Z}|_{[0,1]}}$ denotes the parallel translation along $\gamma_{Z}|_{[0,1]}$ and 
$(\Delta'_V)_{Z}$ and $(\Delta'_H)_{Z}$ are given by 
\begin{align*}
(\Delta'_V)_{Z} &:= \{\lambda \in (\Delta'_V)_+ \mid 0 < \lambda(Z) < \pi \},\\
(\Delta'_H)_{Z} &:= \left\{ \lambda \in (\Delta'_H)_+ \mid -\frac{\pi}{2} < \lambda(Z) < \frac{\pi}{2} \right\}.
\end{align*}

\item[{\rm (ii)}] $H \cdot p\ (p:= \exp_{eK}(Z))$ is a principal orbit if and only if $Z \in C'$ holds.

\item[{\rm (iii)}] Let $\rho_p^{\perp}$ be the normal slice representation of $G \curvearrowright G/K$ at $p$.  Then 
$$T_p^{\perp}(H \cdot p)=\rho_p^{\perp}(H_p) \cdot P_{\gamma_{Z}|_{[0,1]}}(\mathfrak{b})$$
holds, where $H_p$ denotes the isotropy subgroup of $H$ at $p$.  Set 
$\xi := \rho_p^{\perp}(h)(P_{\gamma_{Z}|_{[0,1]}}(\vv)) \in T_p^{\perp}(H \cdot p)\ 
(h \in H_p, \vv \in \mathfrak{b})$.  
If dim$H \cdot p \geq 1$, then the shape operator $A^p_{\xi}$ of the submanifold $H\cdot p$ for $\xi$ satisfies 
\begin{align*}
A^p_{\xi}|_{dh_p(P_{\gamma_{Z}|_{[0,1]}}(\mathfrak{p}_{\lambda} \cap \mathfrak{q}))} &= -\frac{\lambda(\vv)}{{\rm tan}\lambda (Z)}\ {\rm id}\quad (\lambda \in (\Delta'_V)_{Z}),\\
A^p_{\xi}|_{dh_p(P_{\gamma_{Z}|_{[0,1]}}(\mathfrak{p}_{\lambda} \cap \mathfrak{h}))} &= \lambda(\vv)\ {\rm tan}\lambda (Z)\ {\rm id}\quad (\lambda \in (\Delta'_H)_{Z}).
\end{align*}

\end{description}}
\vspace{15pt}

In particular, in the case where $H \cdot p$ ($p:= \exp_{eK}(Z)$) is a principal orbit, that is, $Z\in C'$, 
we have the following fact.\\
\\
\textbf{Proposition 1.2.2.}\ \ {\sl 
Let $p:= \exp_{eK}(Z) \in \exp_{eK}(C')\ (Z \in C')$.  
\begin{description}[labelwidth=24pt] 
\item[{\rm (i)}]  The tangent space $T_p(H\cdot p)$ is decomposed orthogonally as follows:
\begin{align*}
T_p(H \cdot p) 
&= (T_p(H \cdot p) \cap P_{\gamma_{Z}|_{[0,1]}}(T_{eK}F)) \oplus (T_p(H \cdot p) \cap T_p(F^{\perp}))\\
&= P_{\gamma_{Z}|_{[0,1]}} \left(\mathfrak{z}_{\mathfrak{p} \cap \mathfrak{h}}(\mathfrak{b}) \oplus 
\left( \underset{\lambda \in (\Delta'_{H})_+}{\oplus}(\mathfrak{p}_{\lambda} \cap \mathfrak{h})\right)
\oplus 
\left( \underset{\lambda \in (\Delta'_{V})_+}{\oplus}(\mathfrak{p}_{\lambda} \cap \mathfrak{q})\right)\right).
\end{align*}

\item[{\rm (ii)}] $T_p^{\perp}(H \cdot p)=P_{\gamma_{Z}|_{[0,1]}}(\mathfrak{b})$ holds. 
The shape operator $A^p_{P_{\gamma_{Z}|_{[0,1]}}(\vv)}$ of the submanifold $H\cdot p$ for 
$P_{\gamma_{Z}|_{[0,1]}}(\vv)$ satisfies 
\begin{align*}
A^p_{P_{\gamma_{Z}|_{[0,1]}}(\vv)}|_{(P_{\gamma_{Z}|_{[0,1]}}(\mathfrak{p}_{\lambda} \cap \mathfrak{q}))} &= -\frac{\lambda(\vv)}{{\rm tan}\lambda (Z)}\ {\rm id}\quad (\lambda \in (\Delta'_V)_+),\\
A^p_{P_{\gamma_{Z}|_{[0,1]}}(\vv)}|_{(P_{\gamma_{Z}|_{[0,1]}}(\mathfrak{p}_{\lambda} \cap \mathfrak{h}))} &= \lambda(\vv)\ {\rm tan}\lambda (Z)\ {\rm id}\quad (\lambda \in (\Delta'_H)_+).
\end{align*}

\end{description}}
\vspace{20pt}

%

\section*{$\S2.$ The mean curvature flow for orbits of Hermann actions}
\addcontentsline{toc}{part}{$\S2.$ The mean curvature flow for principal orbits of Hermann actions}

In this section, we use the notations in Subsection 1.2.  
Let $H\curvearrowright G/K$ be a commuting Hermann action, and $\mathfrak b$ and $C'$ as in Subsection 1.2.  
Denote by $\vH^{Z}$ the mean curvature vector field of the orbit $H(\exp_{eK}(Z))$ ($Z\in\overline{C'}$).  
Define a vector field $\vX$ on $C'$ by 
$$\vX_{Z}:=(\vH^Z)_{\exp_{eK}(Z)}\quad\,\,(Z\in\mathcal C'),$$
where we note that $T^{\perp}_{\exp_{eK}(Z)}(H(\exp_{eK}(Z)))$ is equal to $T_{\exp_{eK}(Z)}\exp_{eK}(\mathfrak b)$ 
and that $T_{\exp_{eK}(Z)}\exp_{eK}(\mathfrak b)$ is identified with $\mathfrak b$.  
Also, for each simplex $\sigma$ costructing the boundary $\partial C'$ (which is a simplicial complex), define a vector field 
$\vX^{\sigma}$ on $\sigma$ by 
$$(\vX^{\sigma})_{Z}:=(\vH^Z)_{\exp_{eK}(Z)}\quad\,\,(Z\in\sigma),$$
where we note that $(\vH^Z)_{\exp_{eK}(Z)}$ is tangent to $\sigma$.  
According to Lemma 3.1 of \cite{Koi}, the mean curvature flow $\{M^Z_t\}_{t\in[0,T)}$ starting from the principal orbit 
$H(\exp_{eK}(Z))$ ($Z\in C'$) is given by 
$$M^Z_t:=H(\exp_{eK}(c^Z(t)))\quad\,\,(t\in[0,T)),$$
where $c^Z$ is the integral curve of $\vX$ starting from $Z$.  
Thus the mean curvature flow starting from the principal orbits of the Hermann action $H\curvearrowright G/K$ is dominated 
by $\vX$.  Hence we suffice to describe $\vX$ explicitly to investigate the behaviour of the mean curvature flow starting from 
the principal orbits.  
Also, the mean curvature flow $\{M^{Z'}_t\}_{t\in[0,T)}$ starting from the singular orbit $H(\exp_{eK}(Z')$ 
($Z'\in\sigma$) is given by 
$$M^{Z'}_t:=H(\exp_{eK}(c^{Z'}(t)))\quad\,\,(t\in[0,T)),$$
where $c^{Z'}$ is the integral curve of $\vX^{\sigma}$ starting from $Z'$.  
Thus the mean curvature flow starting from the singular orbits through $\exp_{eK}(\sigma)$ of 
the Hermann action $H\curvearrowright G/K$ is dominated by $\vX^{\sigma}$.  
Hence we suffice to describe $\vX^{\sigma}$ explicitly to investigate the behaviour of the mean curvature flow starting from 
the singular orbits through $\exp_{eK}(\sigma)$.  

Set $m_{\lambda}^V:={\rm dim}(\mathfrak p_{\lambda}\cap\mathfrak q)$ ($\lambda\in(\triangle'_V)_+$) and 
$m_{\lambda}^H:={\rm dim}(\mathfrak p_{\lambda}\cap\mathfrak h)$ ($\lambda\in(\triangle'_H)_+$).  
According to Propositions 1.2.1 and 1.2.2, we have 
$$\vX_Z=-\sum_{\lambda\in(\triangle'_V)_+}m_{\lambda}^V\,\cot\lambda(Z)\vv_{\lambda}
+\sum_{\lambda\in(\triangle'_H)_+}m_{\lambda}^H\,\tan\lambda(Z)\vv_{\lambda}\leqno{(2.1)}$$
and 
$$(\vX^{\sigma})_Z=-\sum_{\lambda\in(\triangle'_V)_{\sigma}}m_{\lambda}^V\,\cot\lambda(Z)\vv_{\lambda}
+\sum_{\lambda\in(\triangle'_H)_{\sigma}}m_{\lambda}^H\,\tan\lambda(Z)\vv_{\lambda},\leqno{(2.2)}$$
where $\vv_{\lambda}$ is the element of $\mathfrak b$ defined by $\lambda(\bullet)=\langle\vv_{\lambda},\bullet\rangle$ and 
$(\triangle'_V)_{\sigma}$ and $(\triangle'_H)_{\sigma}$ is given by 
$(\triangle'_V)_{\sigma}:=(\triangle'_V)_{Z}$ ($Z\in\sigma$) and $(\triangle'_H)_{Z}:=(\triangle'_H)_{Z}$ ($Z\in\sigma$).  
Here we note that $(\triangle'_V)_{Z}$ and $(\triangle'_H)_{Z}$ are independent of the choice of $Z\in\sigma$.  

We consider the case where $G/K$ is of rank two and the commuting Hermann action $H\curvearrowright G/K$ is of cohomogeneity two.  
Such commuting Hermann actions $H\curvearrowright G/K$ are classified as in Table 2.1 and 
$(\triangle'_+,(\Delta'V)_+,(\Delta'_H)_+)$ are given as in Table 2.2.  Also, according to $(2.1)$, $\vX$ is given as in Table 2.3.  
These facts were shown in \cite{Koi}.  
In Table 1, $H^{\ast}\curvearrowright G^{\ast}/K$ implies 
the dual action of $H\curvearrowright G/K$ and $L^{\ast}/H\cap K$ is the dual 
of $L/H\cap K$, where $L:={\rm Fix}(\theta\circ\tau)$.  
In Table 2.2, $\{\alpha,\beta,\alpha+\beta\}$ implies a positive root system 
of the root system of (${\mathfrak a}_2$)-type 
($\vv_{\alpha}=(2,0),\vv_{\beta}=(-1,\sqrt3)$), 
$\{\alpha,\beta,\alpha+\beta, 2\alpha+\beta\}$ implies a positive root system 
of the root system of (${\mathfrak b}_2$)(=(${\mathfrak c}_2$))-type 
($\vv_{\alpha}=(1,0),\vv_{\beta}=(-1,1)$), 
$\{\alpha,\beta,\alpha+\beta,2\alpha+\beta,2\alpha,2\alpha+2\beta\}$ implies a positive root system of the root system of 
($\mathfrak{bc}_{2}$)-type ($\vv_{\alpha}=(1,0), \vv_{\beta} = (-1,1)$) 
and 
$\{\alpha,\beta,\alpha+\beta,\alpha+2\beta,\alpha+3\beta,2\alpha+3\beta\}$ 
implies a positive root system 
of the root system of (${\mathfrak g}_2$)-type 
($\vv_{\alpha}=(2\sqrt3,0),\vv_{\beta}=(-\sqrt3,1)$).  
In Table $2.1\sim 2.3$, $\rho_i$ ($i=1,\cdots,16$) imply 
automorphisms of $G$ and $(\cdot)^2$ implies the product Lie group 
$(\cdot)\times(\cdot)$ of a Lie group $(\cdot)$.  In Table 2.2, 
$\displaystyle{\mathop{\alpha}_{(m)}}$ and so on imply that 
the multiplicity of $\alpha$ is equal to $m$.  
Here we identify $\mathfrak b$ with $\mathbb R^2$ under an orthonormal basis $(\ve_1,\ve_2)$ of $\mathfrak b$ 
(i.e., $\displaystyle{(x_1,x_2) (\in\mathbb R^2)\,\,\mathop{=}_{\rm ident.}\,\,x_1\ve_1+x_2\ve_2 (\in\mathfrak b)}$).  

\vspace{0.5truecm}


\begin{longtable}{|c|c|c|}
    \caption*{TABLE 2.1.}\label{tab:long}\\

\hline
{\scriptsize$H\curvearrowright G/K$} & {\scriptsize$H^{\ast}\curvearrowright 
G^{\ast}/K$} & {\scriptsize$L^{\ast}/H\cap K$}\\
\hline
{\scriptsize$\rho_1(SO(3))\curvearrowright SU(3)/SO(3)$} & 
{\scriptsize$SO_0(1,2)\curvearrowright SL(3,{\Bbb R})/SO(3)$} & 
{\scriptsize$(SL(2,{\Bbb R})/SO(2))\times{\Bbb R}$}\\
\hline
{\scriptsize$SO(6)\curvearrowright SU(6)/Sp(3)$} & 
{\scriptsize$SO^{\ast}(6)\curvearrowright SU^{\ast}(6)/Sp(3)$} & 
{\scriptsize$SL(3,{\Bbb C})/SU(3)$}\\
\hline
{\scriptsize$\rho_2(Sp(3))\curvearrowright SU(6)/Sp(3)$} & 
{\scriptsize$Sp(1,2)\curvearrowright SU^{\ast}(6)/Sp(3)$} & 
{\scriptsize$(SU^{\ast}(4)/Sp(2))\times U(1)$}\\
\hline
{\scriptsize$SO(q+2)\curvearrowright$} 
& {\scriptsize$SO_0(2,q)\curvearrowright$} & 
{\scriptsize$SO_0(2,q)/SO(2)\times SO(q)$}\\
{\scriptsize$SU(q+2)/S(U(2)\times U(q))$} 
& {\scriptsize$SU(2,q)/S(U(2)\times U(q))$} & 
{\scriptsize}\\
\hline
{\scriptsize$S(U(j+1)\times U(q-j+1))\curvearrowright$} 
& {\scriptsize$S(U(1,j)\times U(1,q-j))\curvearrowright$} & 
{\scriptsize$(SU(1,j)/S(U(1)\times U(j)))\times$}\\
{\scriptsize$SU(q+2)/S(U(2)\times U(q))$} 
& {\scriptsize$SU(2,q)/S(U(2)\times U(q))$} & 
{\scriptsize$(SU(1,q-j)/S(U(1)\times U(q-j)))$}\\
\hline
{\scriptsize$SO(j+1)\times SO(q-j+1)\curvearrowright$} 
& {\scriptsize$SO(1,j)\times SO(1,q-j)\curvearrowright$} & 
{\scriptsize$(SO_0(1,j)/SO(j))\times$}\\
{\scriptsize$SO(q+2)/SO(2)\times SO(q)$} 
& {\scriptsize$SO(2,q)/SO(2)\times SO(q)$} & 
{\scriptsize$(SO_0(1,q-j)/SO(q-j))$}\\
\hline
{\scriptsize$SO(4)\times SO(4)\curvearrowright$} 
& {\scriptsize$SO^{\ast}(4)\times SO^{\ast}(4)\curvearrowright$} & 
{\scriptsize$SU(2,2)/S(U(2)\times U(2))$}\\
{\scriptsize$SO(8)/U(4)$} 
& {\scriptsize$SO^{\ast}(8)/U(4)$} & 
{\scriptsize}\\
\hline
{\scriptsize$\rho_3(SO(4)\times SO(4))\curvearrowright$} 
& {\scriptsize$SO(4,{\Bbb C})\curvearrowright SO^{\ast}(8)/U(4)$} & 
{\scriptsize$SO(4,{\Bbb C})/SO(4)$}\\
{\scriptsize$SO(8)/U(4)$} 
& {\scriptsize} & {\scriptsize}\\
\hline
{\scriptsize$\rho_4(U(4))\curvearrowright SO(8)/U(4)$} 
& {\scriptsize$U(2,2)\curvearrowright SO^{\ast}(8)/U(4)$} & 
{\scriptsize$(SO^{\ast}(4)/U(2))\times(SO^{\ast}(4)/U(2))$}\\
\hline
{\scriptsize$SO(4)\times SO(6)\curvearrowright$} 
& {\scriptsize$SO^{\ast}(4)\times SO^{\ast}(6)\curvearrowright$} & 
{\scriptsize$SU(2,3)/S(U(2)\times U(3))$}\\
{\scriptsize$SO(10)/U(5)$} 
& {\scriptsize$SO^{\ast}(10)/U(5)$} & 
{\scriptsize}\\
\hline
{\scriptsize$SO(5)\times SO(5)\curvearrowright$} 
& {\scriptsize$SO^(5,{\Bbb C})\curvearrowright 
SO^{\ast}(10)/U(5)$} & 
{\scriptsize$SO(5,{\Bbb C})/SO(5)$}\\
{\scriptsize$SO(10)/U(5)$} 
& {\scriptsize} & {\scriptsize}\\
\hline
{\scriptsize$\rho_5(U(5))\curvearrowright SO(10)/U(5)$} 
& {\scriptsize$U(2,3)\curvearrowright 
SO^{\ast}(10)/U(5)$} & 
{\scriptsize$(SO^{\ast}(4)/U(2))\times(SO^{\ast}(6)/U(3))$}\\
\hline
{\scriptsize$SO(2)^2\times SO(3)^2\curvearrowright$} 
& {\scriptsize$SO(2,{\Bbb C})\times SO(3,{\Bbb C})\curvearrowright$} & 
{\scriptsize$SO_0(2,3)/SO(2)\times SO(3)$}\\
{\scriptsize$(SO(5)\times SO(5))/SO(5)$} 
& {\scriptsize$SO(5,{\Bbb C})/SO(5)$} & 
{\scriptsize}\\
\hline
{\scriptsize$\rho_6(SO(5))\curvearrowright$} 
& {\scriptsize$SO_0(2,3)\curvearrowright SO(5,{\Bbb C})/SO(5)$} & 
{\scriptsize$(SO(2,{\Bbb C})/SO(2))$}\\
{\scriptsize$(SO(5)\times SO(5))/SO(5)$} 
& {\scriptsize} & 
{\scriptsize$\times(SO(3,{\Bbb C})/SO(3))$}\\
\hline
{\scriptsize$\rho_7(U(2))\curvearrowright Sp(2)/U(2)$} 
& {\scriptsize$U(1,1)\curvearrowright Sp(2,{\Bbb R})/U(2)$} & 
{\scriptsize$(Sp(1,{\Bbb R})/U(1))$}\\
{\scriptsize} & {\scriptsize} & {\scriptsize$\times(Sp(1,{\Bbb R})/U(1))$}\\
\hline
{\scriptsize$SU(q+2)\curvearrowright$} 
& {\scriptsize$SU(2,q)\curvearrowright$} & 
{\scriptsize$SU(2,q)/S(U(2)\times U(q))$}\\
{\scriptsize$Sp(q+2)/Sp(2)\times Sp(q)$} 
& {\scriptsize$Sp(2,q)/Sp(2)\times Sp(q)$} & 
{\scriptsize}\\
\hline
{\scriptsize$H\curvearrowright G/K$} & {\scriptsize$H^{\ast}\curvearrowright 
G^{\ast}/K$} & {\scriptsize$L^{\ast}/H\cap K$}\\
\hline
{\scriptsize$U(4)\curvearrowright$} 
& {\scriptsize$U^{\ast}(4)\curvearrowright$} & 
{\scriptsize$Sp(2,{\Bbb C})/Sp(2)$}\\
{\scriptsize$Sp(4)/Sp(2)\times Sp(2)$} 
& {\scriptsize$Sp(2,2)/Sp(2)\times Sp(2)$} & 
{\scriptsize}\\
\hline

\newpage 

\hline
{\scriptsize$Sp(j+1)\times Sp(q-j+1)\curvearrowright$} 
& {\scriptsize$Sp(1,j)\times Sp(1,q-j)\curvearrowright$} & 
{\scriptsize$(Sp(1,j)/Sp(1)\times Sp(j))\times$}\\
{\scriptsize$Sp(q+2)/Sp(2)\times Sp(q)$} 
& {\scriptsize$Sp(2,q)/Sp(2)\times Sp(q)$} & 
{\scriptsize$(Sp(1,q-j)/Sp(1)\times Sp(q-j))$}\\
\hline
{\scriptsize$SU(2)^2\cdot SO(2)^2\curvearrowright$} 
& {\scriptsize$SL(2,{\Bbb C})\cdot SO(2,{\Bbb C})\curvearrowright$} & 
{\scriptsize$Sp(2,{\Bbb R})/U(2)$}\\
{\scriptsize$(Sp(2)\times Sp(2))/Sp(2)$} 
& {\scriptsize$Sp(2,{\Bbb C})/Sp(2)$} & {\scriptsize}\\
\hline
{\scriptsize$\rho_8(Sp(2))\curvearrowright$} 
& {\scriptsize$Sp(2,{\Bbb R})\curvearrowright Sp(2,{\Bbb C})/Sp(2)$} & 
{\scriptsize$(SL(2,{\Bbb C})/SU(2))$}\\
{\scriptsize$(Sp(2)\times Sp(2))/Sp(2)$} 
& {\scriptsize} & {\scriptsize$\times(SO(2,{\Bbb C})/SO(2))$}\\
\hline
{\scriptsize$\rho_9(Sp(2))\curvearrowright$} 
& {\scriptsize$Sp(1,1)\curvearrowright Sp(2,{\Bbb C})/Sp(2)$} & 
{\scriptsize$(Sp(1,{\Bbb C})/Sp(1))$}\\
{\scriptsize$(Sp(2)\times Sp(2))/Sp(2)$} 
& {\scriptsize} & {\scriptsize$\times(Sp(1,{\Bbb C})/Sp(1))$}\\
\hline
{\scriptsize$Sp(4)\curvearrowright E_6/Spin(10)\cdot U(1)$} 
& {\scriptsize$Sp(2,2)\curvearrowright E_6^{-14}/Spin(10)\cdot U(1)$} & 
{\scriptsize$Sp(2,2)/Sp(2)\times Sp(2)$}\\
\hline
{\scriptsize$SU(6)\cdot SU(2)\curvearrowright$} 
& {\scriptsize$SU(2,4)\cdot SU(2)\curvearrowright$} & 
{\scriptsize$SU(2,4)/S(U(2)\times U(4))$}\\
{\scriptsize$ E_6/Spin(10)\cdot U(1)$} 
& {\scriptsize$E_6^{-14}/Spin(10)\cdot U(1)$} & 
{\scriptsize}\\
\hline
{\scriptsize$\rho_{10}(SU(6)\cdot SU(2))\curvearrowright$} 
& {\scriptsize$SU(1,5)\cdot SL(2,{\Bbb R})\curvearrowright$} & 
{\scriptsize$SO^{\ast}(10)/U(5)$}\\
{\scriptsize$ E_6/Spin(10)\cdot U(1)$} 
& {\scriptsize$E_6^{-14}/Spin(10)\cdot U(1)$} & 
{\scriptsize}\\
\hline
{\scriptsize$\rho_{11}(Spin(10)\cdot U(1))\curvearrowright$} 
& {\scriptsize$SO^{\ast}(10)\cdot U(1)\curvearrowright$} & 
{\scriptsize$(SU(1,5)/S(U(1)\times U(5)))$}\\
{\scriptsize$ E_6/Spin(10)\cdot U(1)$} 
& {\scriptsize$E_6^{-14}/Spin(10)\cdot U(1)$} & 
{\scriptsize$\times(SL(2,{\Bbb R})/SO(2))$}\\
\hline
{\scriptsize$\rho_{12}(Spin(10)\cdot U(1))\curvearrowright$} 
& {\scriptsize$SO_0(2,8)\cdot U(1)\curvearrowright$} & 
{\scriptsize$SO_0(2,8)/SO(2)\times SO(8)$}\\
{\scriptsize$ E_6/Spin(10)\cdot U(1)$} 
& {\scriptsize$E_6^{-14}/Spin(10)\cdot U(1)$} & 
{\scriptsize}\\
\hline
{\scriptsize$Sp(4)\curvearrowright E_6/F_4$} 
& {\scriptsize$Sp(1,3)\curvearrowright E_6^{-26}/F_4$} & 
{\scriptsize$SU^{\ast}(6)/Sp(3)$}\\
\hline
{\scriptsize$\rho_{13}(F_4)\curvearrowright E_6/F_4$} 
& {\scriptsize$F_4^{-20}\curvearrowright E_6^{-26}/F_4$} & 
{\scriptsize$(SO_0(1,9)/SO(9))\times U(1)$}\\
\hline
{\scriptsize$\rho_{14}(SO(4))\curvearrowright G_2/SO(4)$} 
& {\scriptsize$SL(2,{\Bbb R})\times SL(2,{\Bbb R})
\curvearrowright G_2^2/SO(4)$} & 
{\scriptsize$SO(4)/SO(2)\times SO(2)$}\\
\hline
{\scriptsize$\rho_{15}(SO(4))\curvearrowright G_2/SO(4)$} 
& {\scriptsize$\rho^{\ast}_{15}(SO(4))\curvearrowright G_2^2/SO(4)$} & 
{\scriptsize$(SL(2,{\Bbb R})/SO(2))$}\\
{\scriptsize} & {\scriptsize} & 
{\scriptsize$\times(SL(2,{\Bbb R})/SO(2))$}\\
\hline
{\scriptsize$\rho_{16}(G_2)\curvearrowright(G_2\times G_2)/G_2$} 
& {\scriptsize$G_2^2\curvearrowright G_2^{\bf C}/G_2$} & 
{\scriptsize$(SL(2,{\Bbb C})/SU(2))$}\\
{\scriptsize} & {\scriptsize} & 
{\scriptsize$\times(SL(2,{\Bbb C})/SU(2))$}\\
\hline
{\scriptsize$SU(2)^4\curvearrowright(G_2\times G_2)/G_2$} 
& {\scriptsize$SL(2,{\Bbb C})\times SL(2,{\Bbb C})\curvearrowright 
G_2^{\bf C}/G_2$} & 
{\scriptsize$G_2^2/SO(4)$}\\
\hline

\end{longtable}

\begin{longtable}{|c|c|c|c|}
    \caption*{TABLE 2.2.}\label{tab:long}\\ 

\hline
$H \curvearrowright G/K$&$\Delta_{+}=\Delta'_{+}$&$(\Delta'_V)_+$&$(\Delta'_H)_+$\\ \hline

\begin{tabular}{c}
$\rho_{1}(SO(3)) \curvearrowright $\\$SU(3)/SO(3)$
\end{tabular}
&
$\{\underset{(1)}{\alpha},\underset{(1)}{\beta},\underset{(1)}{\alpha+\beta}\}$
&
$\{\underset{(1)}{\alpha}\}$
&
$\{\underset{(1)}{\beta},\underset{(1)}{\alpha+\beta}\}$\\ \hline

\begin{tabular}{c}
$SO(6) \curvearrowright $\\$SU(6)/Sp(3)$
\end{tabular}
&
$\{\underset{(4)}{\alpha},\underset{(4)}{\beta},\underset{(4)}{\alpha+\beta}\}$
&
$\{\underset{(2)}{\alpha},\underset{(2)}{\beta},\underset{(2)}{\alpha+\beta}\}$
&
$\{\underset{(2)}{\alpha},\underset{(2)}{\beta},\underset{(2)}{\alpha+\beta}\}$\\ \hline

\begin{tabular}{c}
$\rho_{2}(Sp(3)) \curvearrowright $\\$SU(6)/Sp(3)$
\end{tabular}
&
$\{\underset{(4)}{\alpha},\underset{(4)}{\beta},\underset{(4)}{\alpha+\beta}\}$
&
$\{\underset{(4)}{\alpha}\}$
&
$\{\underset{(4)}{\beta},\underset{(4)}{\alpha+\beta}\}$\\ \hline

\begin{tabular}{c}
$SO(q+2) \curvearrowright $\\$SU(q+2)/$\\$S(U(2)\times U(q))$\\$(q>2)$
\end{tabular}
&
\begin{tabular}{c}
$\{\underset{(2q-4)}{\alpha},\underset{(2)}{\beta},$\\$\underset{(2q-4)}{\alpha+\beta},\underset{(2)}{2\alpha+\beta},$\\$\underset{(1)}{2\alpha},\underset{(2)}{2\alpha+2\beta}\}$
\end{tabular}
&
\begin{tabular}{c}
$\{\underset{(q-2)}{\alpha},\underset{(1)}{\beta},$\\$\underset{(q-2)}{\alpha+\beta},$\\$\underset{(1)}{2\alpha+\beta}\}$
\end{tabular}
&
\begin{tabular}{c}
$\{\underset{(2q-4)}{\alpha},\underset{(2)}{\beta},$\\$\underset{(2q-4)}{\alpha+\beta},\underset{(2)}{2\alpha+\beta},$\\$\underset{(1)}{2\alpha},\underset{(2)}{2\alpha+2\beta}\}$
\end{tabular}\\ \hline

\begin{tabular}{c}
$S(U(j+1)\times $\\$U(q-j+1))\curvearrowright $\\$SU(q+2)/$\\$S(U(2)\times U(q))$\\$(q>2)$
\end{tabular}
&
\begin{tabular}{c}
$\{\underset{(2q-4)}{\alpha},\underset{(2)}{\beta},$\\$\underset{(2q-4)}{\alpha+\beta},\underset{(2)}{2\alpha+\beta},$\\$\underset{(1)}{2\alpha},\underset{(2)}{2\alpha+2\beta}\}$
\end{tabular}
&
\begin{tabular}{c}
$\{\underset{(2j-2)}{\alpha},$\\$\underset{(2q-2j-2)}{\alpha+\beta},$\\$\underset{(1)}{2\alpha},\underset{(1)}{2\alpha+2\beta}\}$
\end{tabular}
&
\begin{tabular}{c}
$\{\underset{(2q-2j-6)}{\alpha},\underset{(2)}{\beta},$\\$\underset{(2j-2)}{\alpha+\beta},\underset{(2)}{2\alpha+\beta}\}$
\end{tabular}\\ \hline

\begin{tabular}{c}
$S(U(2)\times U(2)) \curvearrowright $\\$SU(4)/$\\$S(U(2)\times U(2))$\\(non-isotropy \\gr. act.)
\end{tabular}
&
\begin{tabular}{c}
$\{\underset{(2)}{\alpha},\underset{(1)}{\beta},\underset{(2)}{\alpha+\beta},$\\$\underset{(1)}{2\alpha+\beta}\}$
\end{tabular}
&
$\{\underset{(1)}{\alpha},\underset{(1)}{\alpha+\beta}\}$
&
\begin{tabular}{c}
$\{\underset{(1)}{\alpha},\underset{(1)}{\beta},\underset{(1)}{\alpha+\beta},$\\$\underset{(1)}{2\alpha+\beta}\}$
\end{tabular}\\ \hline

\begin{tabular}{c}
$SO(j+1)\times $\\$SO(q-j+1)\curvearrowright $\\$SO(q+2)/$\\$SO(2)\times SO(q)$
\end{tabular}
&
\begin{tabular}{c}
$\{\underset{(q-2)}{\alpha},\underset{(1)}{\beta},\underset{(q-2)}{\alpha+\beta},$\\$\underset{(1)}{2\alpha+\beta}\}$
\end{tabular}
&
$\{\underset{(j-1)}{\alpha},\underset{(q-j-1)}{\alpha+\beta}\}$
&
\begin{tabular}{c}
$\{\underset{(q-j-1)}{\alpha},\underset{(1)}{\beta},$\\$\underset{(j-1)}{\alpha+\beta},\underset{(1)}{2\alpha+\beta}\}$
\end{tabular}\\ \hline

\begin{tabular}{c}
$SO(4)\times SO(4) \curvearrowright $\\$SO(8)/U(4)$
\end{tabular}
&
\begin{tabular}{c}
$\{\underset{(4)}{\alpha},\underset{(1)}{\beta},\underset{(4)}{\alpha+\beta},$\\$\underset{(1)}{2\alpha+\beta}\}$
\end{tabular}
&
\begin{tabular}{c}
$\{\underset{(2)}{\alpha},\underset{(1)}{\beta},\underset{(2)}{\alpha+\beta},$\\$\underset{(1)}{2\alpha+\beta}\}$
\end{tabular}
&
$\{\underset{(2)}{\alpha},\underset{(2)}{\alpha+\beta}\}$\\ \hline

\begin{tabular}{c}
$\rho_{3}(SO(4)\times SO(4))$\\$ \curvearrowright $\\$SO(8)/U(4)$
\end{tabular}
&
\begin{tabular}{c}
$\{\underset{(4)}{\alpha},\underset{(1)}{\beta},\underset{(4)}{\alpha+\beta},$\\$\underset{(1)}{2\alpha+\beta}\}$
\end{tabular}
&
$\{\underset{(2)}{\alpha},\underset{(2)}{\alpha+\beta}\}$
&
\begin{tabular}{c}
$\{\underset{(2)}{\alpha},\underset{(1)}{\beta},\underset{(2)}{\alpha+\beta},$\\$\underset{(1)}{2\alpha+\beta}\}$
\end{tabular}\\ \hline

\begin{tabular}{c}
$\rho_{4}(U(4)) \curvearrowright $\\$SO(8)/U(4)$
\end{tabular}
&
\begin{tabular}{c}
$\{\underset{(4)}{\alpha},\underset{(1)}{\beta},\underset{(4)}{\alpha+\beta},$\\$\underset{(1)}{2\alpha+\beta}\}$
\end{tabular}
&
$\{\underset{(1)}{\alpha},\underset{(1)}{\alpha+\beta}\}$
&
\begin{tabular}{c}
$\{\underset{(3)}{\alpha},\underset{(1)}{\beta},\underset{(3)}{\alpha+\beta},$\\$\underset{(1)}{2\alpha+\beta}\}$
\end{tabular}\\ \hline

\begin{tabular}{c}
$SO(4)\times SO(6) \curvearrowright $\\$SO(10)/U(5)$
\end{tabular}
&
\begin{tabular}{c}
$\{\underset{(4)}{\alpha},\underset{(4)}{\beta},\underset{(4)}{\alpha+\beta},$\\$\underset{(4)}{2\alpha+\beta},\underset{(1)}{2\alpha},$\\$\underset{(1)}{2\alpha+2\beta}\}$
\end{tabular}
&
\begin{tabular}{c}
$\{\underset{(2)}{\alpha},\underset{(2)}{\beta},\underset{(2)}{\alpha+\beta},$\\$\underset{(2)}{2\alpha+\beta},\underset{(1)}{2\alpha},$\\$\underset{(1)}{2\alpha+2\beta}\}$
\end{tabular}
&
\begin{tabular}{c}
$\{\underset{(2)}{\alpha},\underset{(2)}{\beta},\underset{(2)}{\alpha+\beta},$\\$\underset{(2)}{2\alpha+\beta}\}$
\end{tabular}
\\ \hline

\begin{tabular}{c}
$SO(5)\times SO(5) \curvearrowright $\\$SO(10)/U(5)$
\end{tabular}
&
\begin{tabular}{c}
$\{\underset{(4)}{\alpha},\underset{(4)}{\beta},\underset{(4)}{\alpha+\beta},$\\$\underset{(4)}{2\alpha+\beta},\underset{(1)}{2\alpha},$\\$\underset{(1)}{2\alpha+2\beta}\}$
\end{tabular}
&
\begin{tabular}{c}
$\{\underset{(2)}{\alpha},\underset{(2)}{\beta},\underset{(2)}{\alpha+\beta},$\\$\underset{(2)}{2\alpha+\beta}\}$
\end{tabular}
&
\begin{tabular}{c}
$\{\underset{(2)}{\alpha},\underset{(2)}{\beta},\underset{(2)}{\alpha+\beta},$\\$\underset{(2)}{2\alpha+\beta},\underset{(1)}{2\alpha},$\\$\underset{(1)}{2\alpha+2\beta}\}$
\end{tabular}
\\ \hline

\begin{tabular}{c}
$\rho_{5}(U(5)) \curvearrowright $\\$SO(10)/U(5)$
\end{tabular}
&
\begin{tabular}{c}
$\{\underset{(4)}{\alpha},\underset{(4)}{\beta},\underset{(4)}{\alpha+\beta},$\\$\underset{(4)}{2\alpha+\beta},\underset{(1)}{2\alpha},$\\$\underset{(1)}{2\alpha+2\beta}\}$
\end{tabular}
&
\begin{tabular}{c}
$\{\underset{(4)}{\alpha},\underset{(1)}{2\alpha},$\\$\underset{(1)}{2\alpha+2\beta}\}$
\end{tabular}
&
\begin{tabular}{c}
$\{\underset{(4)}{\beta},\underset{(4)}{\alpha+\beta},$\\$\underset{(4)}{2\alpha+\beta}\}$
\end{tabular}
\\ \hline

\begin{tabular}{c}
$SO(2)^2\times SO(3)^2 $\\$\curvearrowright $\\$(SO(5)\times SO(5))$\\$/SO(5)$
\end{tabular}
&
\begin{tabular}{c}
$\{\underset{(2)}{\alpha},\underset{(2)}{\beta},\underset{(2)}{\alpha+\beta},$\\$\underset{(2)}{2\alpha+\beta}\}$
\end{tabular}
&
\begin{tabular}{c}
$\{\underset{(1)}{\alpha},\underset{(1)}{\beta},\underset{(1)}{\alpha+\beta},$\\$\underset{(1)}{2\alpha+\beta}\}$
\end{tabular}
&
\begin{tabular}{c}
$\{\underset{(1)}{\alpha},\underset{(1)}{\beta},\underset{(1)}{\alpha+\beta},$\\$\underset{(1)}{2\alpha+\beta}\}$
\end{tabular}
\\ \hline

\begin{tabular}{c}
$\rho_{6}(SO(5)) \curvearrowright $\\$(SO(5)\times SO(5))/$\\$SO(5)$
\end{tabular}
&
\begin{tabular}{c}
$\{\underset{(2)}{\alpha},\underset{(2)}{\beta},\underset{(2)}{\alpha+\beta},$\\$\underset{(2)}{2\alpha+\beta}\}$
\end{tabular}
&
$\{\underset{(2)}{\alpha}\}$
&
\begin{tabular}{c}
$\{\underset{(2)}{\beta},\underset{(2)}{\alpha+\beta},$\\$\underset{(2)}{2\alpha+\beta}\}$
\end{tabular}
\\ \hline

\begin{tabular}{c}
$\rho_{7}(U(2)) \curvearrowright $\\$Sp(2)/U(2)$
\end{tabular}
&
\begin{tabular}{c}
$\{\underset{(1)}{\alpha},\underset{(1)}{\beta},\underset{(1)}{\alpha+\beta},$\\$\underset{(1)}{2\alpha+\beta}\}$
\end{tabular}
&
$\{\underset{(1)}{\alpha},\underset{(1)}{\alpha+\beta}\}$
&
$\{\underset{(1)}{\beta},\underset{(1)}{2\alpha+\beta}\}$

\\ \hline

\begin{tabular}{c}
$SU(q+2) \curvearrowright $\\$Sp(q+2)/$\\$Sp(2)\times Sp(q)$\\$(q>2)$
\end{tabular}
&
\begin{tabular}{c}
$\{\underset{(4q-8)}{\alpha},\underset{(4)}{\beta},$\\$\underset{(4q-8)}{\alpha+\beta},\underset{(4)}{2\alpha+\beta},$\\$\underset{(3)}{2\alpha},\underset{(3)}{2\alpha+2\beta}\}$
\end{tabular}
&
\begin{tabular}{c}
$\{\underset{(2q-4)}{\alpha},\underset{(2)}{\beta},$\\$\underset{(2q-4)}{\alpha+\beta},\underset{(2)}{2\alpha+\beta},$\\$\underset{(1)}{2\alpha},\underset{(1)}{2\alpha+2\beta}\}$
\end{tabular}
&
\begin{tabular}{c}
$\{\underset{(2q-4)}{\alpha},\underset{(2)}{\beta},$\\$\underset{(2q-4)}{\alpha+\beta},\underset{(2)}{2\alpha+\beta},$\\$\underset{(2)}{2\alpha},\underset{(2)}{2\alpha+2\beta}\}$
\end{tabular}
\\ \hline

\begin{tabular}{c}
$SU(4) \curvearrowright $\\$Sp(4)/$\\$Sp(2)\times Sp(2)$
\end{tabular}
&
\begin{tabular}{c}
$\{\underset{(4)}{\alpha},\underset{(3)}{\beta},\underset{(3)}{\alpha+\beta},$\\$\underset{(4)}{2\alpha+\beta}\}$
\end{tabular}
&
\begin{tabular}{c}
$\{\underset{(2)}{\alpha},\underset{(1)}{\beta},\underset{(2)}{\alpha+\beta},$\\$\underset{(1)}{2\alpha+\beta}\}$
\end{tabular}
&
\begin{tabular}{c}
$\{\underset{(2)}{\alpha},\underset{(2)}{\beta},\underset{(1)}{\alpha+\beta},$\\$\underset{(3)}{2\alpha+\beta}\}$
\end{tabular}
\\ \hline

\begin{tabular}{c}
$U(4) \curvearrowright $\\$Sp(4)/$\\$Sp(2)\times Sp(2)$
\end{tabular}
&
\begin{tabular}{c}
$\{\underset{(4)}{\alpha},\underset{(3)}{\beta},\underset{(3)}{\alpha+\beta},$\\$\underset{(4)}{2\alpha+\beta}\}$
\end{tabular}
&
\begin{tabular}{c}
$\{\underset{(2)}{\alpha},\underset{(2)}{\beta},\underset{(2)}{\alpha+\beta},$\\$\underset{(2)}{2\alpha+\beta}\}$
\end{tabular}
&
\begin{tabular}{c}
$\{\underset{(2)}{\alpha},\underset{(1)}{\beta},\underset{(1)}{\alpha+\beta},$\\$\underset{(2)}{2\alpha+\beta}\}$
\end{tabular}
\\ \hline

\begin{tabular}{c}
$Sp(j+1)\times$\\$Sp(q-j+1)\curvearrowright $\\$(Sp(2)\times Sp(2))/$\\$Sp(2)$
\end{tabular}
&
\begin{tabular}{c}
$\{\underset{(4q-8)}{\alpha},\underset{(4)}{\beta},$\\$\underset{(4q-8)}{\alpha+\beta},\underset{(4)}{2\alpha+\beta},$\\$\underset{(3)}{2\alpha},\underset{(3)}{2\alpha+2\beta}\}$
\end{tabular}
&
\begin{tabular}{c}
$\{\underset{(2j-4)}{\alpha},\underset{(3)}{2\alpha},$\\$\underset{(4q-4j-4)}{\alpha+\beta},$\\$\underset{(3)}{2\alpha+2\beta}\}$
\end{tabular}
&
\begin{tabular}{c}
$\{\underset{(4q-4j-4)}{\alpha},\underset{(4)}{\beta},$\\$\underset{(4j-4)}{\alpha+\beta},$\\$\underset{(4)}{2\alpha+\beta}\}$
\end{tabular}
\\ \hline

\begin{tabular}{c}
$Sp(2)\times Sp(2)\curvearrowright $\\$Sp(4)/Sp(2)\times $\\$Sp(2)$
\end{tabular}
&
\begin{tabular}{c}
$\{\underset{(4)}{\alpha},\underset{(3)}{\beta},\underset{(3)}{\alpha+\beta},$\\$\underset{(4)}{2\alpha+\beta}\}$
\end{tabular}
&
$\{\underset{(3)}{\alpha},\underset{(3)}{\alpha+\beta}\}$
&
\begin{tabular}{c}
$\{\underset{(1)}{\alpha},\underset{(3)}{\beta},$\\$\underset{(4)}{2\alpha+\beta}\}$
\end{tabular}
\\ \hline

\begin{tabular}{c}
$SU(2)^2\cdot SO(2)^2 \curvearrowright $\\$(Sp(2)\times Sp(2))/$\\$Sp(2)$
\end{tabular}
&
\begin{tabular}{c}
$\{\underset{(2)}{\alpha},\underset{(2)}{\beta},\underset{(2)}{\alpha+\beta},$\\$\underset{(2)}{2\alpha+\beta}\}$
\end{tabular}
&
\begin{tabular}{c}
$\{\underset{(1)}{\alpha},\underset{(1)}{\beta},\underset{(1)}{\alpha+\beta},$\\$\underset{(1)}{2\alpha+\beta}\}$
\end{tabular}
&
\begin{tabular}{c}
$\{\underset{(1)}{\alpha},\underset{(1)}{\beta},\underset{(1)}{\alpha+\beta},$\\$\underset{(1)}{2\alpha+\beta}\}$
\end{tabular}
\\ \hline

\begin{tabular}{c}
$\rho_{8}(Sp(2)) \curvearrowright $\\$(Sp(2)\times Sp(2))/$\\$Sp(2)$
\end{tabular}
&
\begin{tabular}{c}
$\{\underset{(2)}{\alpha},\underset{(2)}{\beta},\underset{(2)}{\alpha+\beta},$\\$\underset{(2)}{2\alpha+\beta}\}$
\end{tabular}
&
$\{\underset{(2)}{\alpha},\underset{(2)}{\alpha+\beta}\}$
&
$\{\underset{(2)}{\beta},\underset{(2)}{2\alpha+\beta}\}$
\\ \hline

\begin{tabular}{c}
$\rho_{9}(Sp(2)) \curvearrowright $\\$(Sp(2)\times Sp(2))/$\\$Sp(2)$
\end{tabular}
&
\begin{tabular}{c}
$\{\underset{(2)}{\alpha},\underset{(2)}{\beta},\underset{(2)}{\alpha+\beta},$\\$\underset{(2)}{2\alpha+\beta}\}$
\end{tabular}
&
$\{\underset{(2)}{\alpha},\underset{(2)}{\alpha+\beta}\}$
&
$\{\underset{(2)}{\beta},\underset{(2)}{2\alpha+\beta}\}$
\\ \hline

\begin{tabular}{c}
$Sp(4) \curvearrowright $\\$E_6/Spin(10)\cdot U(1)$
\end{tabular}
&
\begin{tabular}{c}
$\{\underset{(8)}{\alpha},\underset{(6)}{\beta},\underset{(9)}{\alpha+\beta},$\\$\underset{(5)}{2\alpha+\beta},\underset{(1)}{2\alpha},$\\$\underset{(1)}{2\alpha+2\beta}\}$
\end{tabular}
&
\begin{tabular}{c}
$\{\underset{(4)}{\alpha},\underset{(3)}{\beta},\underset{(3)}{\alpha+\beta},$\\$\underset{(4)}{2\alpha+\beta},\underset{(1)}{2\alpha},$\\$\underset{(1)}{2\alpha+2\beta}\}$
\end{tabular}
&
\begin{tabular}{c}
$\{\underset{(4)}{\alpha},\underset{(3)}{\beta},\underset{(6)}{\alpha+\beta},$\\$\underset{(1)}{2\alpha+\beta}\}$
\end{tabular}
\\ \hline

\begin{tabular}{c}
$SU(6)\cdot SU(2) \curvearrowright $\\$E_6/Spin(10)\cdot U(1)$
\end{tabular}
&
\begin{tabular}{c}
$\{\underset{(8)}{\alpha},\underset{(6)}{\beta},\underset{(9)}{\alpha+\beta},$\\$\underset{(5)}{2\alpha+\beta},\underset{(1)}{2\alpha},$\\$\underset{(1)}{2\alpha+2\beta}\}$
\end{tabular}
&
\begin{tabular}{c}
$\{\underset{(4)}{\alpha},\underset{(2)}{\beta},\underset{(4)}{\alpha+\beta},$\\$\underset{(2)}{2\alpha+\beta},\underset{(1)}{2\alpha},$\\$\underset{(1)}{2\alpha+2\beta}\}$
\end{tabular}
&
\begin{tabular}{c}
$\{\underset{(4)}{\alpha},\underset{(4)}{\beta},\underset{(5)}{\alpha+\beta},$\\$\underset{(1)}{2\alpha+\beta}\}$
\end{tabular}
\\ \hline

\begin{tabular}{c}
$\rho_{10}(SU(6)\cdot SU(2))$\\$\curvearrowright $\\$E_6/Spin(10)\cdot U(1)$
\end{tabular}
&
\begin{tabular}{c}
$\{\underset{(8)}{\alpha},\underset{(6)}{\beta},\underset{(9)}{\alpha+\beta},$\\$\underset{(5)}{2\alpha+\beta},\underset{(1)}{2\alpha},$\\$\underset{(1)}{2\alpha+2\beta}\}$
\end{tabular}
&
\begin{tabular}{c}
$\{\underset{(4)}{\alpha},\underset{(4)}{\beta},\underset{(4)}{\alpha+\beta},$\\$\underset{(4)}{2\alpha+\beta},\underset{(1)}{2\alpha},$\\$\underset{(1)}{2\alpha+2\beta}\}$
\end{tabular}
&
\begin{tabular}{c}
$\{\underset{(4)}{\alpha},\underset{(2)}{\beta},\underset{(5)}{\alpha+\beta},$\\$\underset{(1)}{2\alpha+\beta}\}$
\end{tabular}
\\ \hline

\begin{tabular}{c}
$\rho_{11}(Spin(10)\cdot $\\$SU(1))\curvearrowright $\\$E_6/Spin(10)\cdot U(1)$
\end{tabular}
&
\begin{tabular}{c}
$\{\underset{(8)}{\alpha},\underset{(6)}{\beta},\underset{(9)}{\alpha+\beta},$\\$\underset{(5)}{2\alpha+\beta},\underset{(1)}{2\alpha},$\\$\underset{(1)}{2\alpha+2\beta}\}$
\end{tabular}
&
\begin{tabular}{c}
$\{\underset{(8)}{\alpha},\underset{(1)}{2\alpha},$\\$\underset{(1)}{2\alpha+2\beta}\}$
\end{tabular}
&
\begin{tabular}{c}
$\{\underset{(6)}{\beta},\underset{(9)}{\alpha+\beta},$\\$\underset{(5)}{2\alpha+\beta}\}$
\end{tabular}
\\ \hline

\begin{tabular}{c}
$\rho_{12}(Spin(10)\cdot $\\$SU(1))\curvearrowright $\\$E_6/Spin(10)\cdot U(1)$
\end{tabular}
&
\begin{tabular}{c}
$\{\underset{(8)}{\alpha},\underset{(6)}{\beta},\underset{(9)}{\alpha+\beta},$\\$\underset{(5)}{2\alpha+\beta},\underset{(1)}{2\alpha},$\\$\underset{(1)}{2\alpha+2\beta}\}$
\end{tabular}
&
\begin{tabular}{c}
$\{\underset{(6)}{\alpha},\underset{(1)}{\beta},\underset{(6)}{\alpha+\beta},$\\$\underset{(1)}{2\alpha+\beta}\}$
\end{tabular}
&
\begin{tabular}{c}
$\{\underset{(2)}{\alpha},\underset{(5)}{\beta},\underset{(3)}{\alpha+\beta},$\\$\underset{(4)}{2\alpha+\beta},\underset{(1)}{2\alpha},$\\$\underset{(1)}{2\alpha+2\beta}\}$
\end{tabular}
\\ \hline

\begin{tabular}{c}
$Sp(4)\curvearrowright E_6/F_4$
\end{tabular}
&
$\{\underset{(8)}{\alpha},\underset{(8)}{\beta},\underset{(8)}{\alpha+\beta}\}$
&
$\{\underset{(4)}{\alpha},\underset{(4)}{\beta},\underset{(4)}{\alpha+\beta}\}$
&
$\{\underset{(4)}{\alpha},\underset{(4)}{\beta},\underset{(4)}{\alpha+\beta}\}$
\\ \hline

\begin{tabular}{c}
$\rho_{13}(F_4)\curvearrowright E_6/F_4$
\end{tabular}
&
$\{\underset{(8)}{\alpha},\underset{(8)}{\beta},\underset{(8)}{\alpha+\beta}\}$
&
$\{\underset{(8)}{\alpha}\}$
&
$\{\underset{(8)}{\beta},\underset{(8)}{\alpha+\beta}\}$
\\ \hline

\begin{tabular}{c}
$\rho_{14}(SO(4))\curvearrowright $\\$ G_2/SO(4)$
\end{tabular}
&
\begin{tabular}{c}
$\{\underset{(1)}{\alpha},\underset{(1)}{\beta},\underset{(1)}{\alpha+\beta},$\\$\underset{(1)}{2\alpha+\beta},\underset{(1)}{3\alpha+\beta},$\\$\underset{(1)}{3\alpha+2\beta}\}$
\end{tabular}
&
$\{\underset{(1)}{\alpha},\underset{(1)}{3\alpha+2\beta}\}$
&
\begin{tabular}{c}
$\{\underset{(1)}{\beta},\underset{(1)}{\alpha+\beta},$\\$\underset{(1)}{2\alpha+\beta},$\\$\underset{(1)}{3\alpha+\beta}\}$
\end{tabular}
\\ \hline

\begin{tabular}{c}
$\rho_{15}(SO(4))\curvearrowright $\\$ G_2/SO(4)$
\end{tabular}
&
\begin{tabular}{c}
$\{\underset{(1)}{\alpha},\underset{(1)}{\beta},\underset{(1)}{\alpha+\beta},$\\$\underset{(1)}{2\alpha+\beta},\underset{(1)}{3\alpha+\beta},$\\$\underset{(1)}{3\alpha+2\beta}\}$
\end{tabular}
&
$\{\underset{(1)}{\alpha},\underset{(1)}{3\alpha+2\beta}\}$
&
\begin{tabular}{c}
$\{\underset{(1)}{\beta},\underset{(1)}{\alpha+\beta},$\\$\underset{(1)}{2\alpha+\beta},$\\$\underset{(1)}{3\alpha+\beta}\}$
\end{tabular}
\\ \hline

\begin{tabular}{c}
$\rho_{16}(G_2)\curvearrowright $\\$ (G_2\times G_2)/G_2$
\end{tabular}
&
\begin{tabular}{c}
$\{\underset{(2)}{\alpha},\underset{(2)}{\beta},\underset{(2)}{\alpha+\beta},$\\$\underset{(2)}{2\alpha+\beta},\underset{(2)}{3\alpha+\beta},$\\$\underset{(2)}{3\alpha+2\beta}\}$
\end{tabular}
&
$\{\underset{(2)}{\alpha},\underset{(2)}{3\alpha+2\beta}\}$
&
\begin{tabular}{c}
$\{\underset{(2)}{\beta},\underset{(2)}{\alpha+\beta},$\\$\underset{(2)}{2\alpha+\beta},$\\$\underset{(2)}{3\alpha+\beta}\}$
\end{tabular}
\\ \hline

\begin{tabular}{c}
$SU(2)^4\curvearrowright $\\$ (G_2\times G_2)/G_2$
\end{tabular}
&
\begin{tabular}{c}
$\{\underset{(2)}{\alpha},\underset{(2)}{\beta},\underset{(2)}{\alpha+\beta},$\\$\underset{(2)}{2\alpha+\beta},\underset{(2)}{3\alpha+\beta},$\\$\underset{(2)}{3\alpha+2\beta}\}$
\end{tabular}
&
\begin{tabular}{c}
$\{\underset{(1)}{\alpha},\underset{(1)}{\beta},\underset{(1)}{\alpha+\beta},$\\$\underset{(1)}{2\alpha+\beta},\underset{(1)}{3\alpha+\beta},$\\$\underset{(1)}{3\alpha+2\beta}\}$
\end{tabular}
&
\begin{tabular}{c}
$\{\underset{(1)}{\alpha},\underset{(1)}{\beta},\underset{(1)}{\alpha+\beta},$\\$\underset{(1)}{2\alpha+\beta},\underset{(1)}{3\alpha+\beta},$\\$\underset{(1)}{3\alpha+2\beta}\}$
\end{tabular}
\\ \hline

\end{longtable}

\newpage


\begin{longtable}{|c|c|}
    \caption*{TABLE 2.3.}\label{tab:long}\\ 

\hline
$H \curvearrowright G/K$&$X\qquad(\widetilde{C})$\\ \hline

\begin{tabular}{c}
$\rho_{1}(SO(3)) \curvearrowright $\\$SU(3)/SO(3)$
\end{tabular}
&
\begin{tabular}{c}
$X_{(x, y)}=($tan$(x+\sqrt{3}y)-2$cot$2x+$tan$(x-\sqrt{3}y),$\\$
\sqrt{3}$tan$(x+\sqrt{3}y)-\sqrt{3}$tan$(x-\sqrt{3}y))$\\$
(\widetilde{C}:0<x, \frac{x}{\sqrt{3}}-\frac{\pi}{2\sqrt{3}}<y<-\frac{x}{\sqrt{3}}+\frac{\pi}{2\sqrt{3}})$
\end{tabular}\\ \hline

\begin{tabular}{c}
$SO(6) \curvearrowright $\\$SU(6)/Sp(3)$
\end{tabular}
&
\begin{tabular}{c}
$X_{(x, y)}=(-4$cot$2x-2$cot$(x-\sqrt{3}y)-2$cot$(x+\sqrt{3}y)$\\$
+4$tan$2x+2$tan$(x-\sqrt{3}y)+2$tan$(x-\sqrt{3}y),$\\$
2\sqrt{3}$cot$(x-\sqrt{3}y)-2\sqrt{3}$cot$(x+\sqrt{3}y)$\\$
-2\sqrt{3}$tan$(x-\sqrt{3}y)+2\sqrt{3}$tan$(x+\sqrt{3}y))$\\$
(\widetilde{C}:0<x, \frac{x}{\sqrt{3}}<y<-\frac{x}{\sqrt{3}}+\frac{\pi}{2\sqrt{3}})$
\end{tabular}\\ \hline


\begin{tabular}{c}
$\rho_{2}(Sp(3)) \curvearrowright $\\$SU(6)/Sp(3)$
\end{tabular}
&
\begin{tabular}{c}
$X_{(x, y)}=(-8$cot$2x+4$tan$(x-\sqrt{3}y)+2$tan$(x-\sqrt{3}y),$\\$
4\sqrt{3}$tan$(x-\sqrt{3}y)-4\sqrt{3}$tan$(x+\sqrt{3}y))$\\$
(\widetilde{C}:0<x, \frac{x}{\sqrt{3}}-\frac{\pi}{2\sqrt{3}}<y<-\frac{x}{\sqrt{3}}+\frac{\pi}{2\sqrt{3}})$
\end{tabular}\\ \hline

\begin{tabular}{c}
$SO(q+2) \curvearrowright $\\$SU(q+2)/$\\$S(U(2)\times U(q))$\\$(q>2)$
\end{tabular}
&
\begin{tabular}{c}
$X_{(x, y)}=(-(q-2)$cot$x-$cot$(x-y)-$cot$(x+y)$\\$
+(q-2)$tan$x-$tan$(x-y)+$tan$(x+y)+$tan$2x,$\\$
$cot$(x-y)-(q-2)$cot$y-$cot$(x+y)$\\$
-$tan$(x-y)+(q-2)$tan$y+$tan$(x+y)+2$tan$2y)$\\$
(\widetilde{C}:0<x<y<\frac{\pi}{4})$
\end{tabular}\\ \hline

\begin{tabular}{c}
$SO(4) \curvearrowright $\\$SU(4)/$\\$S(U(2)\times U(2))$
\end{tabular}
&
\begin{tabular}{c}
$X_{(x, y)}=(-$cot$x+$tan$x-$tan$(x-y)+$tan$(x+y),$\\$
-$cot$y-$tan$(x-y)+$tan$y+$tan$(x+y))$\\$
(\widetilde{C}:0<x<y<\frac{\pi}{4})$
\end{tabular}\\ \hline

\begin{tabular}{c}
$S(U(j+1)\times $\\$U(q-j+1))\curvearrowright $\\$SU(q+2)/$\\$S(U(2)\times U(q))$\\$(q>2)$
\end{tabular}
&
\begin{tabular}{c}
$X_{(x, y)}=(-2(j-1)$cot$x-2$cot$2x+2(q-j-1)$tan$x$\\$
-2$tan$(x-y)+2$tan$(x+y),$\\$
-2(q-j-1)$cot$y-2$cot$2y-2$tan$(x-y)$\\$
+2(j-1)$tan$y+2$tan$(x+y))$\\$
(\widetilde{C}:0<x<y<\frac{\pi}{4})$
\end{tabular}\\ \hline

\begin{tabular}{c}
$S(U(2)\times U(2)) \curvearrowright $\\$SU(4)/$\\$S(U(2)\times U(2))$\\(non-isotropy \\gr. act.)
\end{tabular}
&
\begin{tabular}{c}
$X_{(x, y)}=(-$cot$x+$tan$x+$tan$(x-y)+$tan$(x+y),$\\$
-$cot$y-$tan$(x-y)+$tan$y+$tan$(x+y))$\\$
(\widetilde{C}:0<x, 0<y, x+y<\frac{\pi}{2})$
\end{tabular}\\ \hline

\begin{tabular}{c}
$SO(j+1)\times $\\$SO(q-j+1)\curvearrowright $\\$SO(q+2)/$\\$SO(2)\times SO(q)$
\end{tabular}
&
\begin{tabular}{c}
$X_{(x, y)}=(-(j-1)$cot$x+(q-j-1)$tan$x$\\$
+$tan$(x-y)+$tan$(x+y),$\\$
-(q-j-1)$cot$y-$tan$(x-y)$\\$
+(j-1)$tan$y+$tan$(x+y))$\\$
(\widetilde{C}:0<x, 0<y, x+y<\frac{\pi}{2})$
\end{tabular}\\ \hline

\begin{tabular}{c}
$SO(4)\times SO(4) \curvearrowright $\\$SO(8)/U(4)$
\end{tabular}
&
\begin{tabular}{c}
$X_{(x, y)}=(-2$cot$x-$cot$(x-y)-2$cot$(x-y)+2$tan$x,$\\$
$cot$(x-y)-2$cot$y-2$cot$(x+y)+2$tan$y)$\\$
(\widetilde{C}:0<y<x, x+y<\pi)$
\end{tabular}\\ \hline

\begin{tabular}{c}
$\rho_{3}(SO(4)\times SO(4))$\\$ \curvearrowright $\\$SO(8)/U(4)$
\end{tabular}
&
\begin{tabular}{c}
$X_{(x, y)}=(-2$cot$x+2$tan$x+$tan$(x-y)+$tan$(x+y),$\\$
-2$cot$y-$tan$(x-y)+2$tan$y+$tan$(x+y))$\\$
(\widetilde{C}:0<x, 0<y, x+y<\frac{\pi}{2})$
\end{tabular}\\ \hline

\begin{tabular}{c}
$\rho_{4}(U(4)) \curvearrowright $\\$SO(8)/U(4)$
\end{tabular}
&
\begin{tabular}{c}
$X_{(x, y)}=(-$cot$x+3$tan$x+$tan$(x-y)+$tan$(x+y),$\\$
-$cot$y-$tan$(x-y)+3$tan$y+$tan$(x+y))$\\$
(\widetilde{C}:0<x, 0<y, x+y<\frac{\pi}{2})$
\end{tabular}\\ \hline


\begin{tabular}{c}
$SO(4)\times SO(6) \curvearrowright $\\$SO(10)/U(5)$
\end{tabular}
&
\begin{tabular}{c}
$X_{(x, y)}=(-$cot$x-$cot$(x-y)-$cot$(x+y)$\\$
-$cot$2x+$tan$x+$tan$(x-y)+$tan$(x+y),$\\$
-$cot$(x-y)-$cot$y-$cot$(x+y)$\\$
-$cot$2y-$tan$(x-y)+$tan$y+$tan$(x+y))$\\$
(\widetilde{C}:0<x<y, x+y<\frac{\pi}{2})$
\end{tabular}\\ \hline

\begin{tabular}{c}
$SO(5)\times SO(5) \curvearrowright $\\$SO(10)/U(5)$
\end{tabular}
&
\begin{tabular}{c}
$X_{(x, y)}=(-$cot$x-$cot$(x-y)-$cot$(x+y)$\\$
+$tan$x+$tan$(x-y)+$tan$(x+y)+$tan$2x,$\\$
-$cot$(x-y)-$cot$y-$cot$(x+y)$\\$
-$tan$(x-y)+$tan$y+$tan$(x+y)+$tan$2y)$\\$
(\widetilde{C}:0<x<y<\frac{\pi}{4})$
\end{tabular}\\ \hline

\begin{tabular}{c}
$\rho_{5}(U(5)) \curvearrowright $\\$SO(10)/U(5)$
\end{tabular}
&
\begin{tabular}{c}
$X_{(x, y)}=2(-2$cot$x-$cot$2x$\\$
+2$tan$(x-y)+2$tan$(x+y),$\\$
-$cot$2y-2$tan$(x-y)+2$tan$(x+y)+2$tan$y)$\\$
(\widetilde{C}:0<x, 0<y, x+y<\frac{\pi}{2})$
\end{tabular}\\ \hline

\begin{tabular}{c}
$SO(2)^2\times SO(3)^2 $\\$\curvearrowright $\\$(SO(5)\times SO(5))$\\$/SO(5)$
\end{tabular}
&
\begin{tabular}{c}
$X_{(x, y)}=(-$cot$x-$cot$(x-y)-$cot$(x+y)$\\$
+$tan$x+$tan$(x-y)+$tan$(x+y),$\\$
$cot$(x-y)-$cot$y-2$cot$(x+y)$\\$
-$tan$(x-y)+$tan$y+$tan$(x+y))$\\$
(\widetilde{C}:0<x<y, x+y<\frac{\pi}{2})$
\end{tabular}\\ \hline

\begin{tabular}{c}
$\rho_{6}(SO(5)) \curvearrowright $\\$(SO(5)\times SO(5))/$\\$SO(5)$
\end{tabular}
&
\begin{tabular}{c}
$X_{(x, y)}=2(-$cot$x+$tan$(x-y)+$tan$(x+y),$\\$
-$tan$(x-y)+$tan$y+$tan$(x+y))$\\$
(\widetilde{C}:0<x, x-\frac{\pi}{2}<y, x+y<\frac{\pi}{2})$
\end{tabular}\\ \hline

\begin{tabular}{c}
$\rho_{7}(U(2)) \curvearrowright $\\$Sp(2)/U(2)$
\end{tabular}
&
\begin{tabular}{c}
$X_{(x, y)}=(-$cot$x+$tan$(x-y)+$tan$(x+y),$\\$
-$cot$y-$tan$(x-y)+$tan$(x+y))$\\$
(\widetilde{C}:0<x, 0<y, x+y<\frac{\pi}{2})$
\end{tabular}\\ \hline

\begin{tabular}{c}
$SU(q+2) \curvearrowright $\\$Sp(q+2)/$\\$Sp(2)\times Sp(q)$\\$(q>2)$
\end{tabular}
&
\begin{tabular}{c}
$X_{(x, y)}=(-(2q-4)$cot$x-2$cot$(x-y)-2$cot$(x+y)$\\$
-2$cot$2x+(2q-4)$tan$x+2$tan$(x-y)$\\$
+2$tan$(x+y)+4$tan$2x,$\\$
2$cot$(x-y)-(2q-4)$cot$y-2$cot$(x+y)$\\$
-2$cot$2y-2$tan$(x-y)+(2q-4)$tan$y$\\$
+2$tan$(x+y)+4$tan$2y)$\\$
(\widetilde{C}:0<x<y<\frac{\pi}{4})$
\end{tabular}\\ \hline

\begin{tabular}{c}
$SU(4) \curvearrowright $\\$Sp(4)/$\\$Sp(2)\times Sp(2)$
\end{tabular}
&
\begin{tabular}{c}
$X_{(x, y)}=(-2$cot$x-$cot$(x-y)-$cot$(x+y)$\\$
+2$tan$x+2$tan$(x-y)+3$tan$(x+y),$\\$
$cot$(x-y)-2$cot$y-$cot$(x+y)$\\$
-2$tan$(x-y)+$tan$y+3$tan$(x+y))$\\$
(\widetilde{C}:0<x<y, x+y<\frac{\pi}{2})$
\end{tabular}\\ \hline


\begin{tabular}{c}
$U(4) \curvearrowright $\\$Sp(4)/$\\$Sp(2)\times Sp(2)$
\end{tabular}
&
\begin{tabular}{c}
$X_{(x, y)}=(-2$cot$x-2$cot$(x-y)-2$cot$(x+y)$\\$
+2$tan$x+$tan$(x-y)+2$tan$(x+y),$\\$
2$cot$(x-y)-2$cot$y-2$cot$(x+y)$\\$
-$tan$(x-y)+$tan$y+2$tan$(x+y))$\\$
(\widetilde{C}:0<x<y, x+y<\frac{\pi}{2})$
\end{tabular}\\ \hline

\begin{tabular}{c}
$Sp(j+1)\times$\\$Sp(q-j+1)\curvearrowright $\\$(Sp(2)\times Sp(2))/$\\$Sp(2)$
\end{tabular}
&
\begin{tabular}{c}
$X_{(x, y)}=(-4(j-1)$cot$x-6$cot$2x$\\$
+4(q-j-1)$tan$x+4$tan$(x-y)+4$tan$(x+y),$\\$
-4(q-j-1)$cot$y-6$cot$2y-4$tan$(x-y)$\\$
+4(j-1)$tan$y+4$tan$(x+y))$\\$
(\widetilde{C}:0<x<y, x+y<\frac{\pi}{2})$
\end{tabular}\\ \hline

\begin{tabular}{c}
$Sp(2)\times Sp(2)\curvearrowright $\\$Sp(4)/Sp(2)\times $\\$Sp(2)$
\end{tabular}
&
\begin{tabular}{c}
$X_{(x, y)}=(-3$cot$x+$tan$x$\\$
+3$tan$(x-y)+4$tan$(x+y),$\\$
-3$cot$y-3$tan$(x-y)+4$tan$(x+y))$\\$
(\widetilde{C}:0<x<y, x+y<\frac{\pi}{2})$
\end{tabular}\\ \hline

\begin{tabular}{c}
$SU(2)^2\cdot SO(2)^2 \curvearrowright $\\$(Sp(2)\times Sp(2))/$\\$Sp(2)$
\end{tabular}
&
\begin{tabular}{c}
$X_{(x, y)}=(-$cot$x-$cot$(x-y)-$cot$(x+y)$\\$
+$tan$x+$tan$(x-y)+$tan$(x+y),$\\$
$cot$(x-y)-$cot$y-$cot$(x+y)$\\$
-$tan$(x-y)+$tan$y+$tan$(x+y))$\\$
(\widetilde{C}:0<x<y, x+y<\frac{\pi}{2})$
\end{tabular}\\ \hline

\begin{tabular}{c}
$\rho_{8}(Sp(2)) \curvearrowright $\\$(Sp(2)\times Sp(2))/$\\$Sp(2)$
\end{tabular}
&
\begin{tabular}{c}
$X_{(x, y)}=2(-$cot$x+$tan$(x-y)+$tan$(x+y),$\\$
-$cot$y-$tan$(x-y)+$tan$(x+y))$\\$
(\widetilde{C}:0<x, 0<y, x+y<\frac{\pi}{2})$
\end{tabular}\\ \hline

\begin{tabular}{c}
$\rho_{9}(Sp(2)) \curvearrowright $\\$(Sp(2)\times Sp(2))/$\\$Sp(2)$
\end{tabular}
&
\begin{tabular}{c}
$X_{(x, y)}=2(-$cot$x+$tan$(x-y)+$tan$(x+y),$\\$
-$cot$y-$tan$(x-y)+$tan$(x+y))$\\$
(\widetilde{C}:0<x, 0<y, x+y<\frac{\pi}{2})$
\end{tabular}\\ \hline

\begin{tabular}{c}
$Sp(4) \curvearrowright $\\$E_6/Spin(10)\cdot U(1)$
\end{tabular}
&
\begin{tabular}{c}
$X_{(x, y)}=(-4$cot$x-3$cot$(x-y)-4$cot$(x+y)$\\$
+4$tan$x+3$tan$(x-y)+$tan$(x+y),$\\$
3$cot$(x-y)-3$cot$y-4$cot$(x+y)$\\$
-3$tan$(x-y)+6$tan$y+$tan$(x+y))$\\$
(\widetilde{C}:0<x<y, x+y<\frac{\pi}{2})$
\end{tabular}\\ \hline

\begin{tabular}{c}
$SU(6)\cdot SU(2) \curvearrowright $\\$E_6/Spin(10)\cdot U(1)$
\end{tabular}
&
\begin{tabular}{c}
$X_{(x, y)}=(-4$cot$x-2$cot$(x-y)-2$cot$(x+y)$\\$
-2$cot$2x+4$tan$x+4$tan$(x-y)+3$tan$(x+y),$\\$
2$cot$(x-y)-4$cot$y-2$cot$(x+y)$\\$
-2$cot$2y-4$tan$(x-y)+5$tan$y+3$tan$(x+y))$\\$
(\widetilde{C}:0<x<y, x+y<\frac{\pi}{2})$
\end{tabular}\\ \hline


\begin{tabular}{c}
$\rho_{10}(SU(6)\cdot SU(2))$\\$\curvearrowright $\\$E_6/Spin(10)\cdot U(1)$
\end{tabular}
&
\begin{tabular}{c}
$X_{(x, y)}=(-4$cot$x-4$cot$(x-y)-4$cot$(x+y)$\\$
-2$cot$2x+4$tan$x+2$tan$(x-y)+$tan$(x+y),$\\$
4$cot$(x-y)-4$cot$y-4$cot$(x+y)$\\$
-2$cot$2y-2$tan$(x-y)+5$tan$y+$tan$(x+y))$\\$
(\widetilde{C}:0<x<y, x+y<\frac{\pi}{2})$
\end{tabular}\\ \hline

\begin{tabular}{c}
$\rho_{11}(Spin(10)\cdot $\\$SU(1))\curvearrowright $\\$E_6/Spin(10)\cdot U(1)$
\end{tabular}
&
\begin{tabular}{c}
$X_{(x, y)}=(-8$cot$x-2$cot$2x$\\$
+6$tan$(x-y)+5$tan$(x+y),$\\$
-2$cot$2y-6$tan$(x-y)+9$tan$y+5$tan$(x+y))$\\$
(\widetilde{C}:0<x, 0<y, x+y<\frac{\pi}{2})$
\end{tabular}\\ \hline

\begin{tabular}{c}
$\rho_{12}(Spin(10)\cdot $\\$SU(1))\curvearrowright $\\$E_6/Spin(10)\cdot U(1)$
\end{tabular}
&
\begin{tabular}{c}
$X_{(x, y)}=(-6$cot$x-$cot$(x-y)-$cot$(x+y)$\\$
-2$cot$2x+5$tan$x+4$tan$(x-y)+2$tan$(x+y),$\\$
$cot$(x-y)-6$cot$y-$cot$(x+y)$\\$
-5$tan$(x-y)+3$tan$y+4$tan$(x+y)+2$tan$2y)$\\$
(\widetilde{C}:0<x<y<\frac{\pi}{4})$
\end{tabular}\\ \hline

\begin{tabular}{c}
$Sp(4)\curvearrowright E_6/F_4$
\end{tabular}
&
\begin{tabular}{c}
$X_{(x, y)}=(-8$cot$2x-4$cot$(x-\sqrt{3}y)-4$cot$(x+\sqrt{3}y)$\\$
+8$tan$2x+4$tan$(x-\sqrt{3}y)+4$tan$(x-\sqrt{3}y),$\\$
4\sqrt{3}$cot$(x-\sqrt{3}y)-4\sqrt{3}$cot$(x+\sqrt{3}y)$\\$
-4\sqrt{3}$tan$(x-\sqrt{3}y)+4\sqrt{3}$tan$(x+\sqrt{3}y))$\\$
(\widetilde{C}:0<x, \frac{x}{\sqrt{3}}<y<-\frac{x}{\sqrt{3}}+\frac{\pi}{2\sqrt{3}})$
\end{tabular}\\ \hline

\begin{tabular}{c}
$\rho_{13}(F_4)\curvearrowright E_6/F_4$
\end{tabular}
&
\begin{tabular}{c}
$X_{(x, y)}=(-16$cot$2x+8$tan$(3x-\sqrt{3}y)+8$tan$(x-\sqrt{3}y),$\\$
-8\sqrt{3}$tan$(x-\sqrt{3}y)+8\sqrt{3}$tan$(x+\sqrt{3}y))$\\$
(\widetilde{C}:0<x, \frac{x}{\sqrt{3}}-\frac{\pi}{2\sqrt{3}}<y<-\frac{x}{\sqrt{3}}+\frac{\pi}{2\sqrt{3}})$
\end{tabular}\\ \hline

\begin{tabular}{c}
$\rho_{14}(SO(4))\curvearrowright $\\$ G_2/SO(4)$
\end{tabular}
&
\begin{tabular}{c}
$X_{(x, y)}=(-2$cot$2x-3$tan$(3x-\sqrt{3}y)-$tan$(x-\sqrt{3}y)$\\$
+$tan$(x+\sqrt{3}y)+3$tan$(3x+\sqrt{3}y),$\\$
-2\sqrt{3}$cot$2\sqrt{3}y-\sqrt{3}$tan$(3x-\sqrt{3}y)-\sqrt{3}$tan$(x-\sqrt{3}y)$\\$
-\sqrt{3}$tan$(x+\sqrt{3}y)+\sqrt{3}$tan$(3x+\sqrt{3}y))$\\$
(\widetilde{C}:0<x, 0<y<-\sqrt{3}x+\frac{\pi}{2\sqrt{3}})$
\end{tabular}\\ \hline

\begin{tabular}{c}
$\rho_{15}(SO(4))\curvearrowright $\\$ G_2/SO(4)$
\end{tabular}
&
\begin{tabular}{c}
$X_{(x, y)}=(-2$cot$2x-3$tan$(3x-\sqrt{3}y)-$tan$(x-\sqrt{3}y)$\\$
+$tan$(x+\sqrt{3}y)+3$tan$(3x+\sqrt{3}y),$\\$
-2\sqrt{3}$cot$2\sqrt{3}y-\sqrt{3}$tan$(3x-\sqrt{3}y)-\sqrt{3}$tan$(x-\sqrt{3}y)$\\$
-\sqrt{3}$tan$(x+\sqrt{3}y)+\sqrt{3}$tan$(3x+\sqrt{3}y))$\\$
(\widetilde{C}:0<x, 0<y<-\sqrt{3}x+\frac{\pi}{2\sqrt{3}})$
\end{tabular}\\ \hline

\begin{tabular}{c}
$\rho_{16}(G_2)\curvearrowright $\\$ (G_2\times G_2)/G_2$
\end{tabular}
&
\begin{tabular}{c}
$X_{(x, y)}=2(-2$cot$2x-3$tan$(3x-\sqrt{3}y)-$tan$(x-\sqrt{3}y)$\\$
+$tan$(x+\sqrt{3}y)+3$tan$(3x+\sqrt{3}y),$\\$
-2\sqrt{3}$cot$2\sqrt{3}y-\sqrt{3}$tan$(3x-\sqrt{3}y)-\sqrt{3}$tan$(x-\sqrt{3}y)$\\$
-\sqrt{3}$tan$(x+\sqrt{3}y)+\sqrt{3}$tan$(3x+\sqrt{3}y))$\\$
(\widetilde{C}:0<x, 0<y<-\sqrt{3}x+\frac{\pi}{2\sqrt{3}})$
\end{tabular}\\ \hline


\begin{tabular}{c}
$SU(2)^4\curvearrowright $\\$ (G_2\times G_2)/G_2$
\end{tabular}
&\begin{tabular}{c}
$X_{(x, y)}=(-2$cot$2x-3$cot$(3x-\sqrt{3}y)-$cot$(x-\sqrt{3}y)$\\$
-$cot$(x+\sqrt{3}y)-3$cot$(3x+\sqrt{3}y)+2$tan$2x$\\$
+3$tan$(3x-\sqrt{3}y)+$tan$(x-\sqrt{3}y)$\\$
+$tan$(x+\sqrt{3}y)+3$tan$(3x+\sqrt{3}y),$\\$
\sqrt{3}$cot$(3x-\sqrt{3}y)-\sqrt{3}$cot$(x-\sqrt{3}y)-\sqrt{3}$cot$(x+\sqrt{3}y)$\\$
-\sqrt{3}$cot$(3x+\sqrt{3}y)-2\sqrt{3}$cot$2\sqrt{3}y-\sqrt{3}$tan$(3x-\sqrt{3}y)$\\$
-\sqrt{3}$tan$(x-\sqrt{3}y)+\sqrt{3}$tan$(x+\sqrt{3}y)$\\$
+\sqrt{3}$tan$(3x+\sqrt{3}y)+2\sqrt{3}$tan$2\sqrt{3}y)$\\$
(\widetilde{C}:0<x, \sqrt{3}x<y<\frac{\pi}{4\sqrt{3}})$
\end{tabular}\\ \hline


\end{longtable}

\vspace{0.5truecm}


\section*{$\S3.$ The minimal orbits of Hermann actions of cohomogeneity two}
\addcontentsline{toc}{part}{$\S3.$ The minimal orbits of Hermann actions of cohomogeneity two}

In this section, we use the notations in Subsection 1.2 and Section 2.  
Let $H\curvearrowright G/K$ be a commuting Hermann action of cohomogeneity two on an irreducible rank two Riemannian symmetric 
space of compact type.  It is known that the only minimal principal $H$-orbit $H(\exp_{eK}(Z_0))$ ($Z_0\in C'$) and 
the only minimal singular $H$-orbit $H(\exp_{eK}(Z^{\sigma}))$ ($Z^{\sigma}\in\sigma$) through the edge $\exp_{eK}(\sigma)$ 
for each edge $\sigma$ of $\partial C'$, where the ``minimality'' means that the orbit is a minimal submanifold in $G/K$.  
Let $Z_0=x_1^0\ve_1+x_2^0\ve_2$ and $Z^{\sigma}=x_1^{\sigma}\ve_1+x_2^{\sigma}\ve_2$, where $(\ve_1,\ve_2)$ is an orthonormal basis 
of $\mathfrak b$ as stated in Section 2.  
It is clear that $\vX_{Z_0}={\bf 0}$ and $(\vX^{\sigma})_{Z^{\sigma}}={\bf 0}$ holds.  Hence, from the descriptions of $\vX$ 
in Table 2.3, we can calculate $(x^0_1,x^0_2)$ by using Mathematicae.  Also, we can describe $\vX^{\sigma}$'s explicitly 
by using the description $(2.2)$ of $\vX^{\sigma}$ and, from the explicit descriptions of $\vX^{\sigma}$, we can calculate 
$(x^{\sigma}_1,x^{\sigma}_2)$ by using Mathematicae.  
In the sequel, we abbreviate $X_{x_1e_1+x_2e_2}$ as $X_{(x_1,x_2)}$.  

For example, we consider the case of the Hermann action $\rho_{1}(SO(3)) \curvearrowright SU(3)/SO(3)$.  
Then, from the description of $\vX$ in Table 2.3, we obtain $(x^0_1,x^0_2)=(\frac{\pi}{6},0)$.  
According to Table 2.2, we have $\Delta_{+} = \{\underset{(1)}{\alpha},\underset{(1)}{\beta},\underset{(1)}{\alpha+\beta}\}$ 
($\mathfrak a_2$-type), $(\Delta_V)_+ = \{\underset{(1)}{\alpha}\}$, $(\Delta_H)_+ = \{\underset{(1)}{\beta},\underset{(1)}{\alpha+\beta}\}$ ($\alpha = (2,0), \beta = (-1,\sqrt{3})$).  
Hence $\sigma$ is equal to an open interval of one of $\alpha^{-1}(0),\,\beta^{-1}(-\frac{\pi}{2})$ and 
$(\alpha+\beta)^{-1}(\frac{\pi}{2})$.  
In the case of $\sigma\subset\alpha^{-1}(0)$, we have 
$$X_{(x_1,x_2)}^{\sigma}
=(0, 2\sqrt{3}\tan \sqrt{3}x_2).$$
Hence we obtain $(x_1^{\sigma}, x_2^{\sigma})=(0,0)$.  

In the case of $\sigma\subset\beta^{-1}(-\frac{\pi}{2})$, we have 
$$X_{(x_1,x_2)}^{\sigma}
=(-3\cot 2\sqrt{3}x_2, -\sqrt{3}\cot 2\sqrt{3}x_2).$$
Hence we obtain $(x_1^{\sigma},x_2^{\sigma}) = (\frac{\pi}{4},-\frac{\pi}{4\sqrt{3}})$.  

In the case of $\sigma\subset(\alpha + \beta)^{-1}(\frac{\pi}{2})$, we have 
$$X_{(x_1,x_2)}^{\sigma}=(3\cot 2\sqrt{3}x_2, -\sqrt{3}\cot 2\sqrt{3}x_2).$$
Hence we obtain $(x_1^{\sigma}, x_2^{\sigma}) = (\frac{\pi}{4},\frac{\pi}{4\sqrt{3}})$.  

In other cases, we can calculate $(x_1^0,x_2^0)$ and $(x_1^{\sigma},x_2^{\sigma})$'s by using Mathematicae.  
As its result, we obtain the following result.  

\vspace{0.5truecm}

\noindent
{\bf Theorem 3.1.}\ \ {\sl Let $H\curvearrowright G/K$ be a commuting Hermann action of cohomogeneity two on an irreducible rank two 
Riemannian symmetric space.  The only minimal principal $H$-orbit is given by $H(\exp_{eK}(x_1^0\ve_1+x_2^0\ve_2))$ and 
the only minimal singular $H$-orbit through $\exp_{eK}(\sigma)$ ($\sigma\,:\,$ one of three edges of $\partial C'$) is given by 
$H(\exp_{eK}(x_1^{\sigma}\ve_1+x_2^{\sigma}\ve_2))$, where $(x_1^0,x_2^0)$ and $(x_1^{\sigma},x_2^{\sigma})$ are as in Table 3.1.  
Also, the above vector field $\vX$ on $C'$ is illustrated as in Table 3.1.}


\newpage

\begin{longtable}{|c|c|c|}
    \caption*{TABLE 3.1.}\label{tab:long}\\ 

\hline
$H \curvearrowright G/K$&$\vX,\,(x_1^0, x_2^0),\,(x_1^{\sigma},x_2^{\sigma})$\\ \hline

\begin{tabular}{c}
$\rho_{1}(SO(3)) \curvearrowright SU(3)/SO(3)$
\end{tabular}
&
\begin{tabular}{c}\\
$\includegraphics[height=50mm]{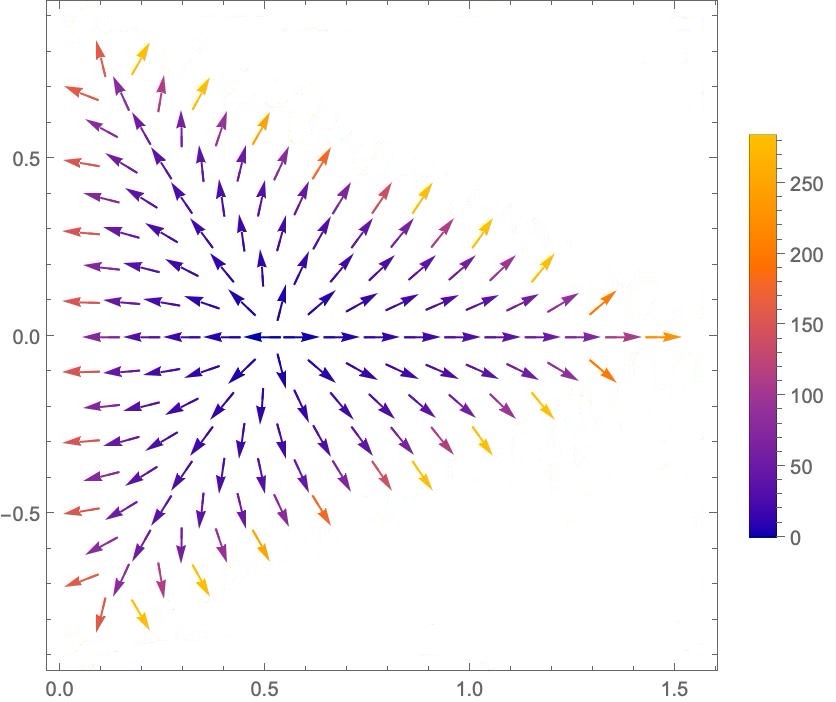}$\\
$(x_1^0,x_2^0) = (\frac{\pi}{6}, 0)$\\
$(x_1^{\sigma}, x_2^{\sigma}) = (0,0), (\frac{\pi}{4},-\frac{\pi}{4\sqrt{3}}), (\frac{\pi}{4},\frac{\pi}{4\sqrt{3}})$
\end{tabular}\\ \hline

\begin{tabular}{c}
$SO(6) \curvearrowright SU(6)/Sp(3)$
\end{tabular}
&
\begin{tabular}{c}\\
$\includegraphics[height=50mm]{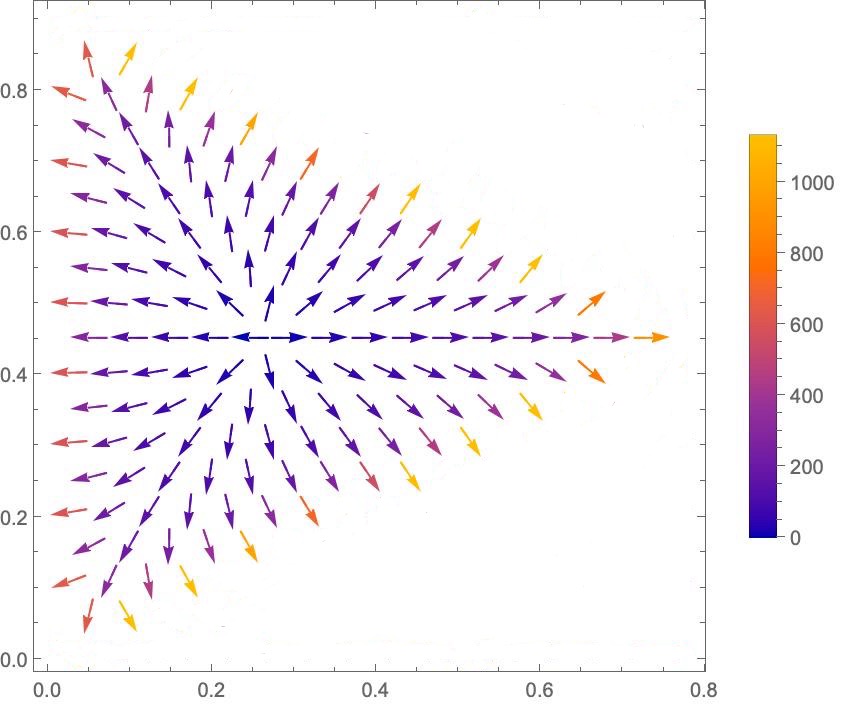}$\\
$(x_1^0,x_2^0) =$ (0.261799, 0.45345)\\
$(x_1^{\sigma},x_2^{\sigma}) = (0,\frac{\pi}{4\sqrt{3}}), (\frac{\pi}{8},\frac{\pi}{8\sqrt{3}}), (\frac{3\pi}{8},\frac{\pi}{8\sqrt{3}})$
\end{tabular}\\ \hline

\begin{tabular}{c}
$\rho_{2}(Sp(3)) \curvearrowright SU(6)/Sp(3)$
\end{tabular}
&
\begin{tabular}{c}\\
$\includegraphics[height=50mm]{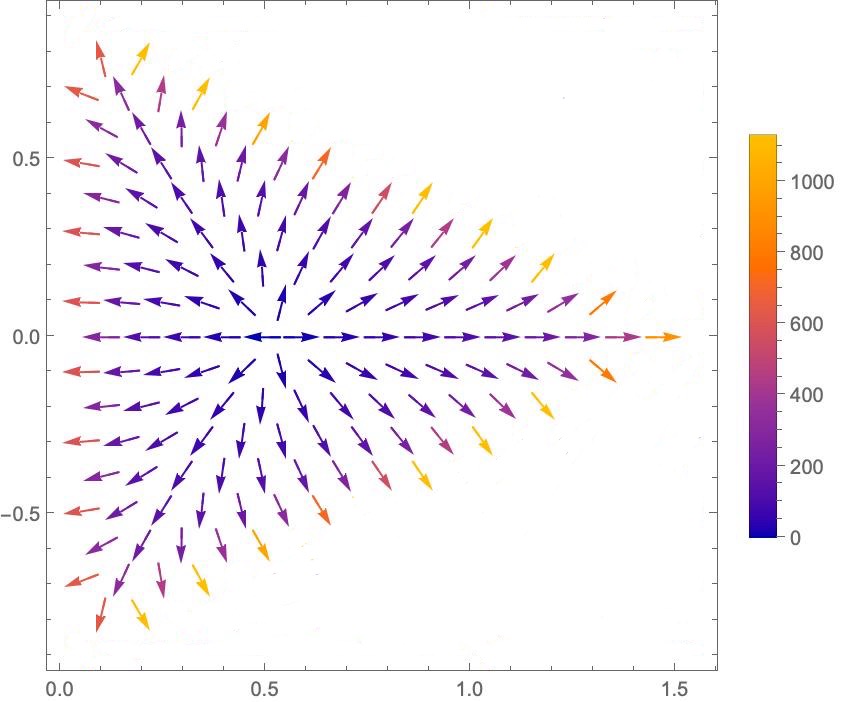}$\\
$(x_1^0,x_2^0) = (\frac{\pi}{6}, 0)$\\
$(x_1^{\sigma},x_2^{\sigma}) = (0,0), (\frac{\pi}{4},-\frac{\pi}{4\sqrt{3}}), (\frac{\pi}{4},\frac{\pi}{4\sqrt{3}})$
\end{tabular}\\ \hline


\hline
\begin{tabular}{c}
$SO(q+2) \curvearrowright $\\$SU(q+2)/S(U(2)\times U(q))$\\$(q=3)$
\end{tabular}
&
\begin{tabular}{c}\\
$\includegraphics[height=50mm]{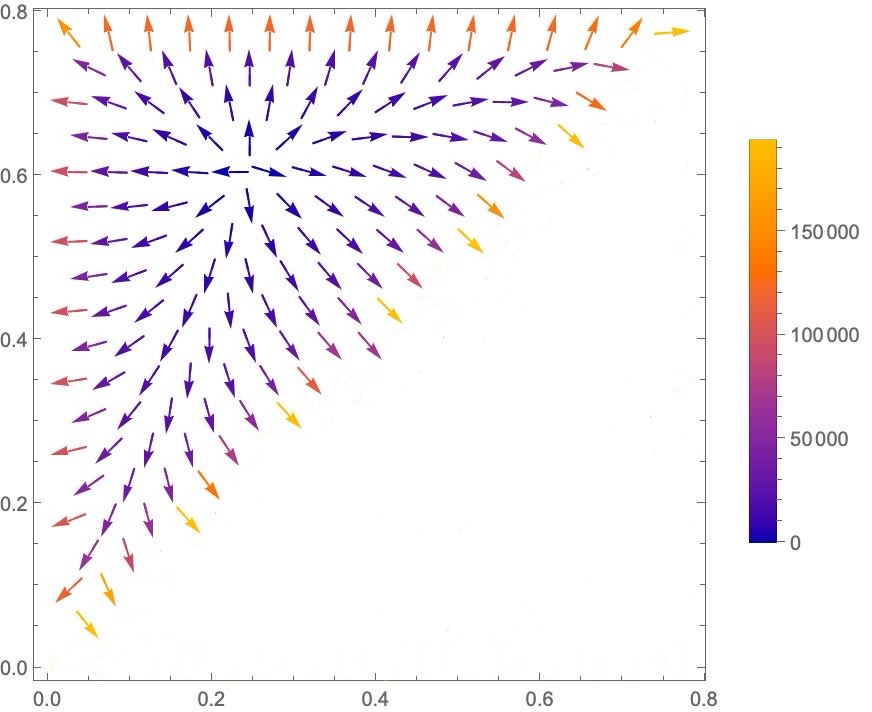}$\\
$(x_1^0,x_2^0) =$ (0.242863, 0.608349)\\
$(x_1^{\sigma},x_2^{\sigma}) = (0,\frac{\pi}{4}), (0,0), (\frac{\pi}{4},\frac{\pi}{4})$
\end{tabular}\\ \hline

\begin{tabular}{c}
$SO(4) \curvearrowright $\\$SU(4)/S(U(2)\times U(2))$
\end{tabular}
&
\begin{tabular}{c}\\
$\includegraphics[height=50mm]{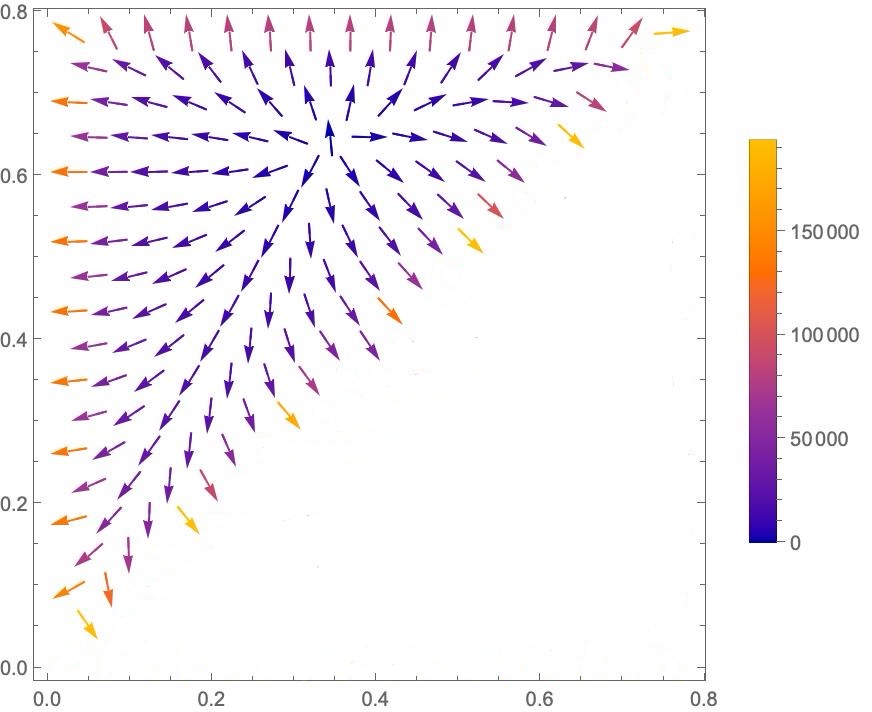}$\\
$(x_1^0,x_2^0) =$ (0.343408, 0.639725)\\
$(x_1^{\sigma},x_2^{\sigma}) =  (0,\arctan(\frac{1}{3})),$\\
$(0.477658,0.477658), (0.33312,\frac{\pi}{4})$
\end{tabular}\\ \hline

\begin{tabular}{c}
$S(U(j+1)\times U(q-j+1))\curvearrowright $\\$SU(q+2)/S(U(2)\times U(q))$\\$(q=3,j=2)$
\end{tabular}
&
\begin{tabular}{c}\\
$\includegraphics[height=50mm]{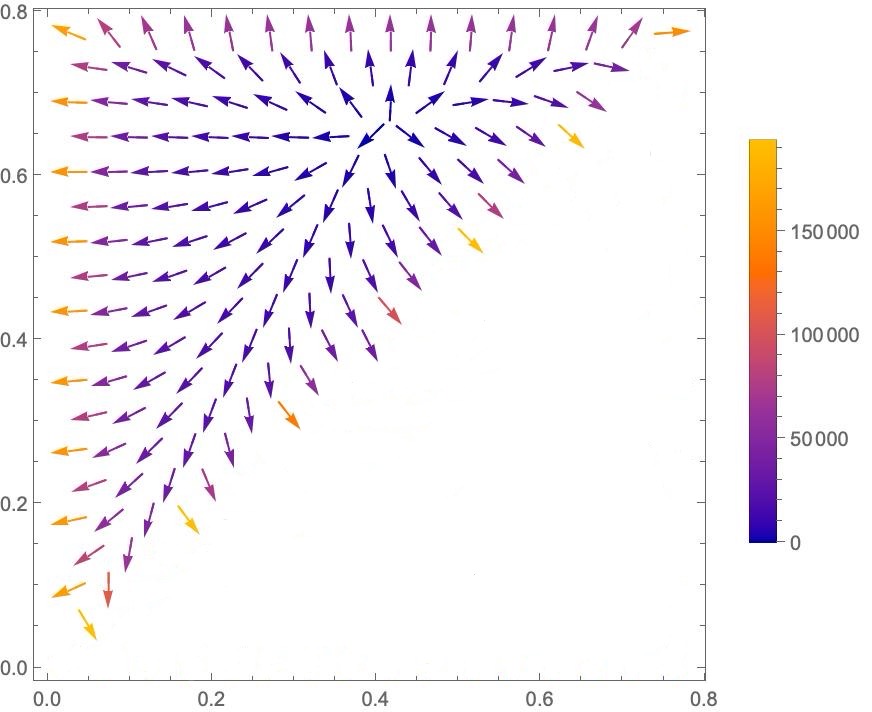}$\\
$(x_1^0,x_2^0) =$ (0.40878, 0.660012)\\
$(x_1^{\sigma},x_2^{\sigma}) = (0,0.31416), (0.560791,0.560791),$\\$ (0.428528,\frac{\pi}{4}))$
\end{tabular}\\ \hline


\hline
\begin{tabular}{c}
$S(U(2)\times U(2)) \curvearrowright $\\$SU(4)/S(U(2)\times U(2))$\\(non-isotropy gr. act.)
\end{tabular}
&
\begin{tabular}{c}\\
$\includegraphics[height=50mm]{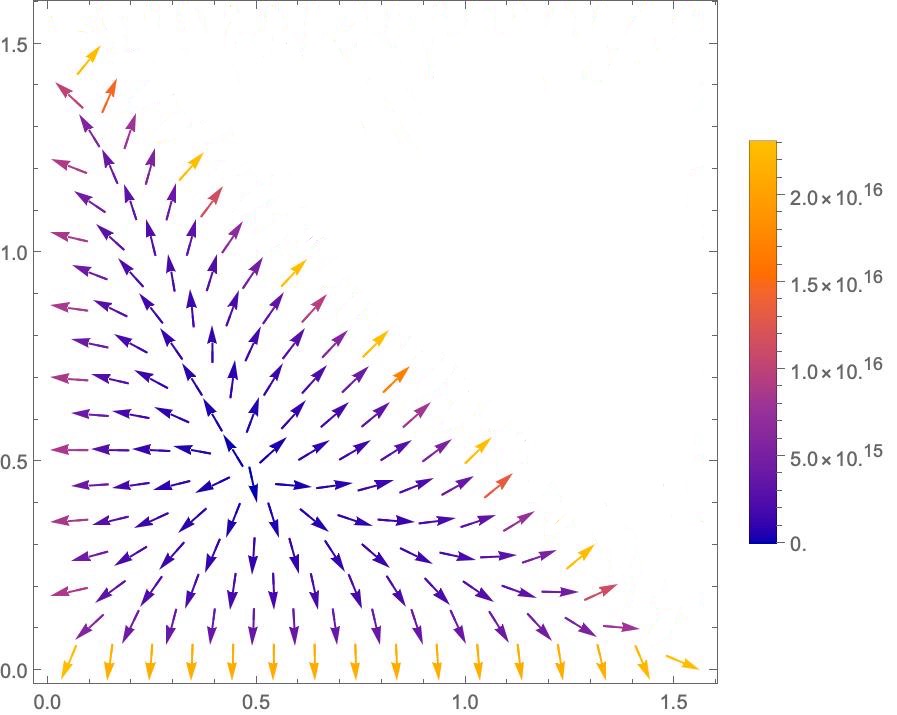}$\\
$(x_1^0,x_2^0) =$ (0.477658, 0.477658)\\
$(x_1^{\sigma},x_2^{\sigma}) = (0,\arctan(\frac{1}{3})), (\arctan(\frac{1}{3}),0), (\frac{\pi}{4},\frac{\pi}{4})$
\end{tabular}\\ \hline

\begin{tabular}{c}
$SO(j+1)\times SO(q-j+1)\curvearrowright $\\$SO(q+2)/SO(2)\times SO(q)$
\end{tabular}
&
\begin{tabular}{c}\\
$\includegraphics[height=50mm]{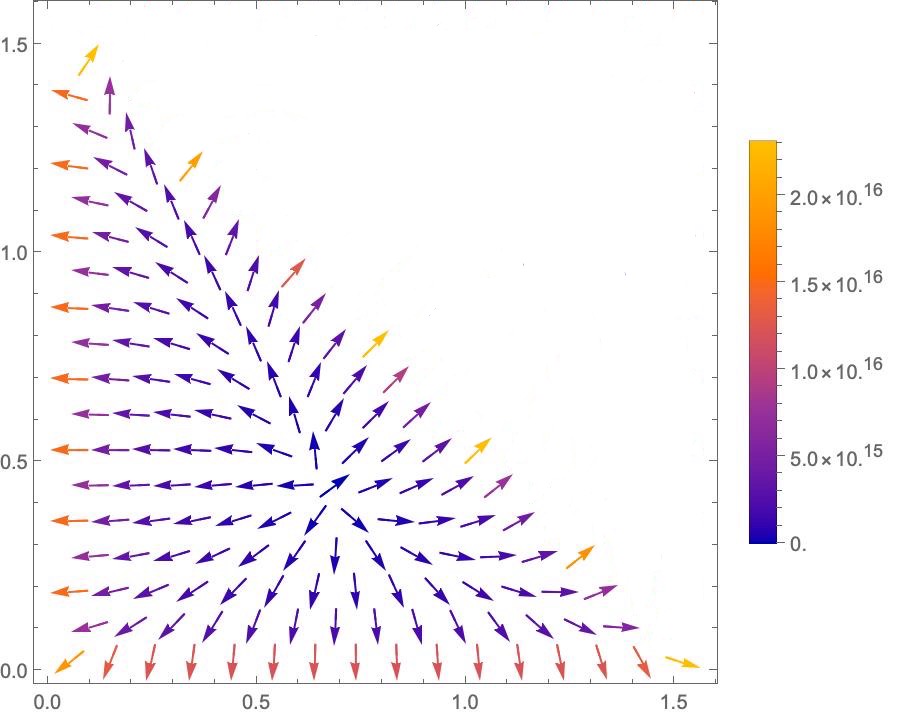}$\\
$(x_1^0,x_2^0) =$ (0.669504, 0.430285)\\
$(x_1^{\sigma},x_2^{\sigma}) = (0,\arctan(\frac{1}{3})), (\arctan(\frac{1}{3}),0), (\frac{\pi}{4},\frac{\pi}{4})$
\end{tabular}\\ \hline

\begin{tabular}{c}
$SO(4)\times SO(4)\curvearrowright $\\$SO(8)/U(4)$
\end{tabular}
&
\begin{tabular}{c}\\
$\includegraphics[height=50mm]{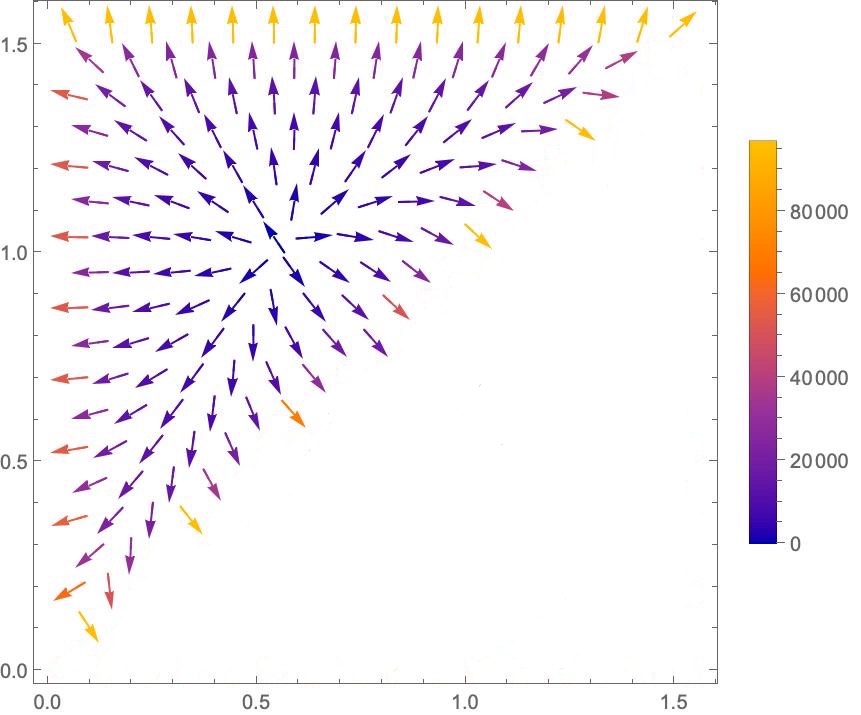}$\\
$(x_1^0,x_2^0) =$ (1.11297, 0.567154)\\
$(x_1^{\sigma},x_2^{\sigma}) = (0,0.85707),(\frac{\pi}{4},\frac{\pi}{4}), (1.00685,\frac{\pi}{2})$
\end{tabular}\\ \hline

%

 \hline
\begin{tabular}{c}
$\rho_{3}(SO(4)\times SO(4)) \curvearrowright $\\$SO(8)/U(4)$
\end{tabular}
&
\begin{tabular}{c}\\
$\includegraphics[height=50mm]{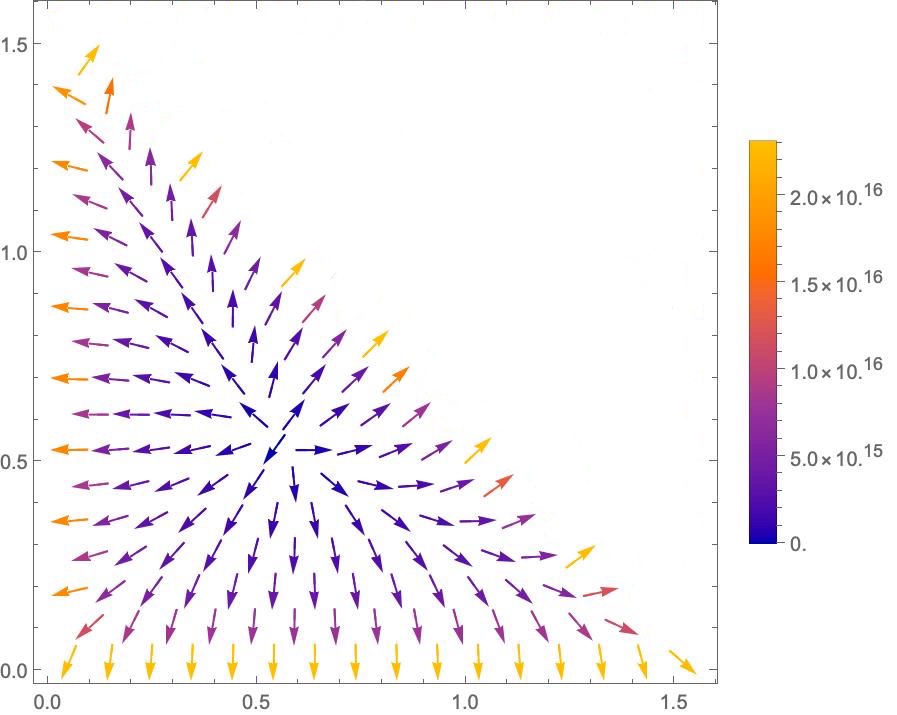}$\\
$(x_1^0,x_2^0) =$ (0.553574, 0.553574)\\
$(x_1^{\sigma},x_2^{\sigma}) = (0,1.00685),(0.61548,0),(\frac{\pi}{4},\frac{\pi}{4})$
\end{tabular}\\ \hline

\begin{tabular}{c}
$\rho_{4}(U(4))\curvearrowright SO(8)/U(4)$
\end{tabular}
&
\begin{tabular}{c}\\
$\includegraphics[height=50mm]{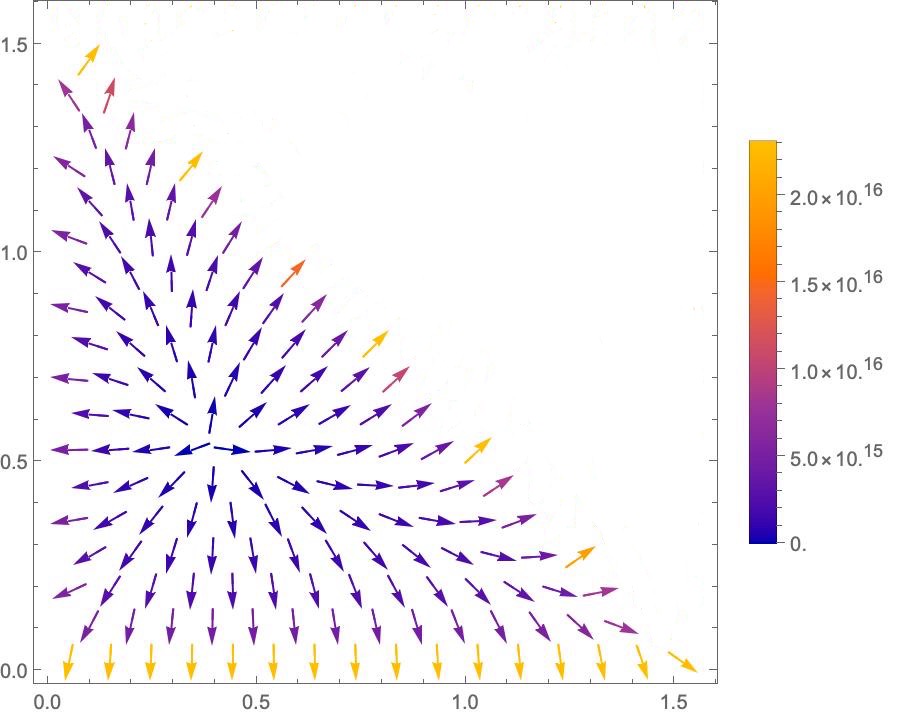}$\\
$(x_1^0,x_2^0) =$ (0.390322, 0.542954)\\
$(x_1^{\sigma},x_2^{\sigma}) = (0,0.42053),(0.42053,0),(\frac{\pi}{4},\frac{\pi}{4})$
\end{tabular}\\ \hline

\begin{tabular}{c}
$SO(4)\times SO(6)\curvearrowright $\\$SO(10)/U(5)$
\end{tabular}
&
\begin{tabular}{c}\\
$\includegraphics[height=50mm]{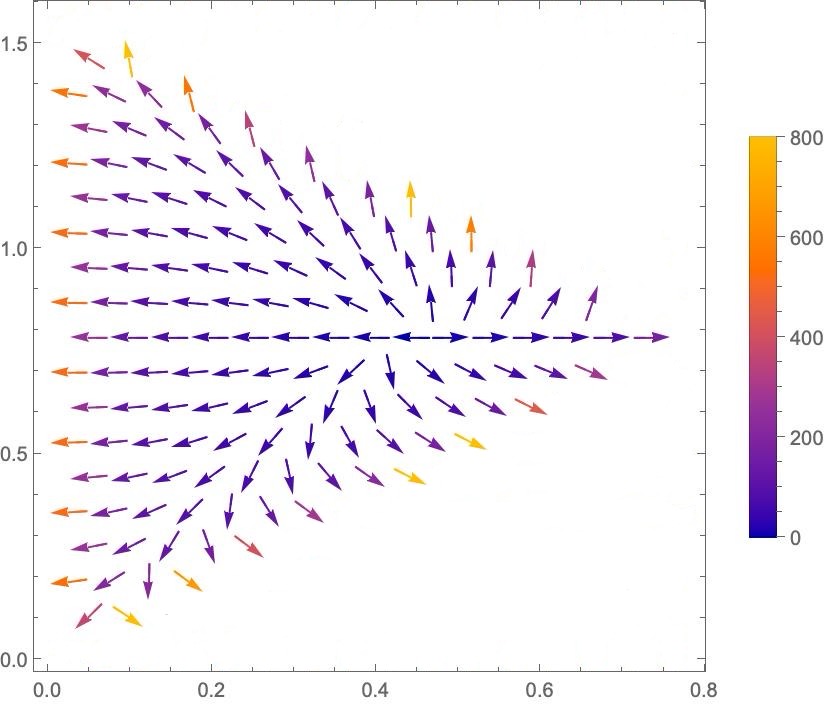}$\\
$(x_1^0,x_2^0) =$ (0.443039, 0.785398)\\
$(x_1^{\sigma},x_2^{\sigma}) = (0,\frac{\pi}{4}), (0,0), (\frac{\pi}{4},\frac{\pi}{4})$
\end{tabular}\\ \hline


\hline
\begin{tabular}{c}
$SO(5)\times SO(5)\curvearrowright $\\$SO(10)/U(5)$
\end{tabular}
&
\begin{tabular}{c}\\
$\includegraphics[height=50mm]{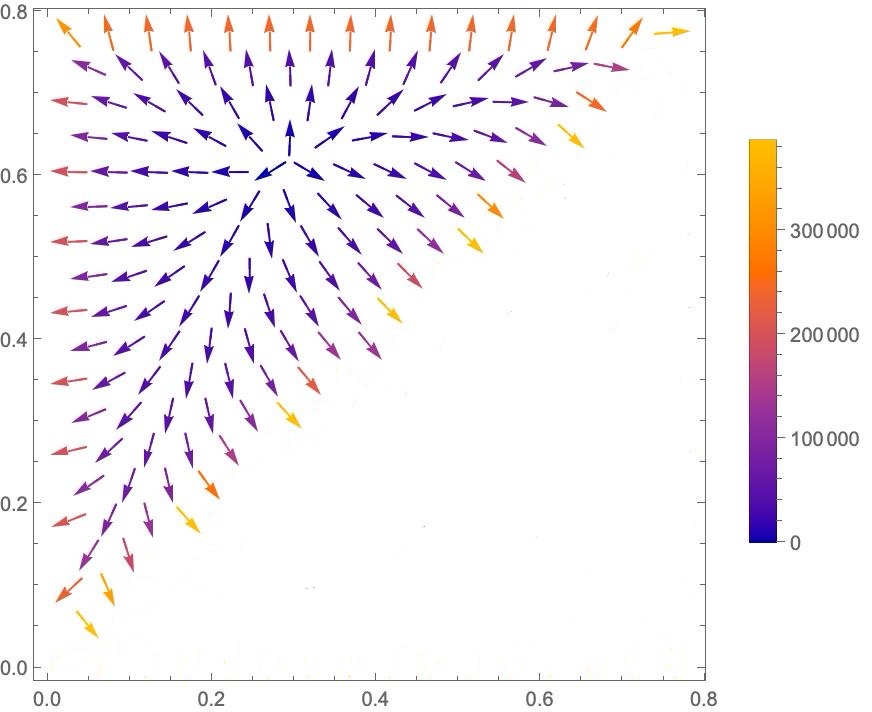}$\\
$(x_1^0,x_2^0) =$ (0.28557, 0.615128)\\
$(x_1^{\sigma},x_2^{\sigma}) = (0,0.5916), (0.44304,0.44304), (0.5916,\frac{\pi}{4})$
\end{tabular}\\ \hline

\begin{tabular}{c}
$\rho_{5}(U(5)) \curvearrowright SO(10)/U(5)$
\end{tabular}
&
\begin{tabular}{c}\\
$\includegraphics[height=50mm]{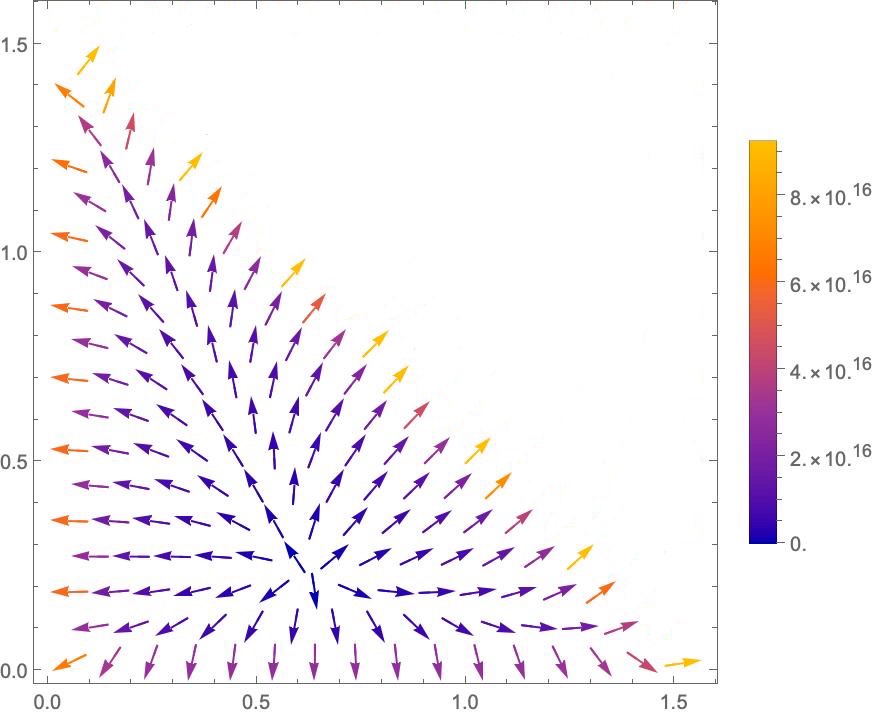}$\\
$(x_1^0,x_2^0) =$ (0.622334, 0.234738)\\
$(x_1^{\sigma},x_2^{\sigma}) = (0,0.36137),(0.64052,0),(0.54453,1.02627)$
\end{tabular}\\ \hline

\begin{tabular}{c}
$SO(2)^2\times SO(3)^2\curvearrowright $\\$(SO(5)\times SO(5))/SO(5)$
\end{tabular}
&
\begin{tabular}{c}\\
$\includegraphics[height=50mm]{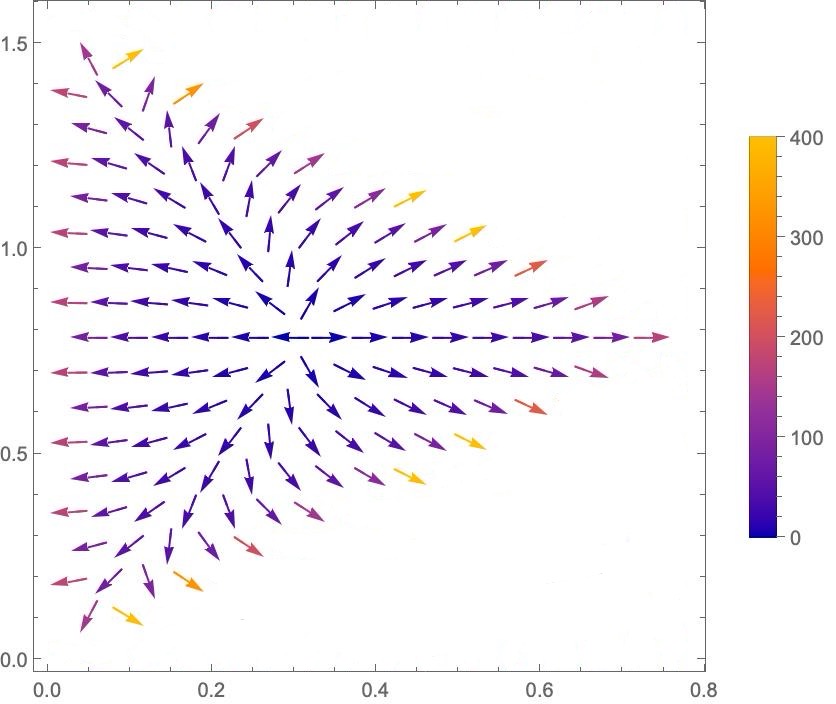}$\\
$(x_1^0,x_2^0) =$ (0.30774, 0.785398)\\
$(x_1^{\sigma},x_2^{\sigma}) = (0,\frac{\pi}{4}), (0,0), (0,\frac{\pi}{2})$
\end{tabular}\\ \hline


\hline
\begin{tabular}{c}
$\rho_{6}(SO(5))\curvearrowright $\\$(SO(5)\times SO(5))/SO(5)$
\end{tabular}
&
\begin{tabular}{c}\\
$\includegraphics[height=50mm]{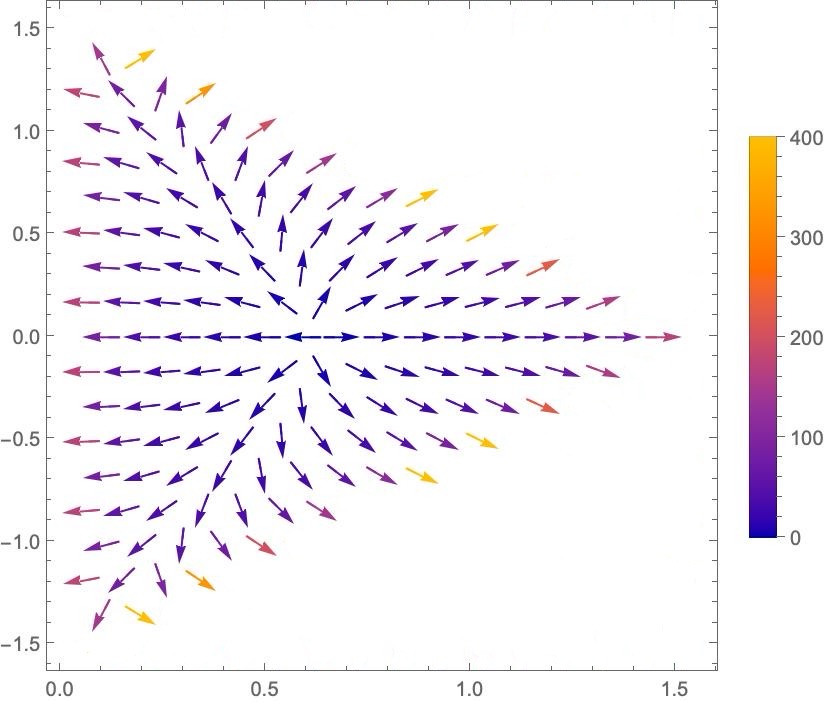}$\\
$(x_1^0,x_2^0) =$ (0.61548, 0)\\
$(x_1^{\sigma},x_2^{\sigma}) = (0,\frac{\pi}{2}), (0,0), (0,-\frac{\pi}{2})$
\end{tabular}\\ \hline

\begin{tabular}{c}
$\rho_{7}(U(2)) \curvearrowright Sp(2)/U(2)$
\end{tabular}
&
\begin{tabular}{c}\\
$\includegraphics[height=50mm]{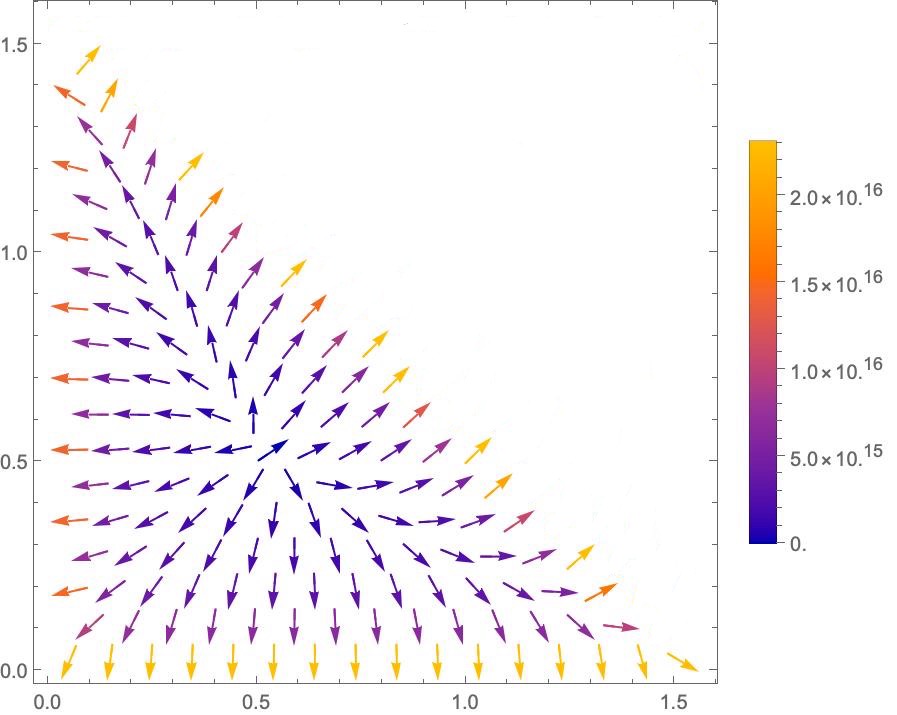}$\\
$(x_1^0,x_2^0) = (\frac{\pi}{6},\frac{\pi}{6})$\\
$(x_1^{\sigma},x_2^{\sigma}) = (0,\frac{\pi}{6}), (\frac{\pi}{6},0), (\frac{\pi}{4},\frac{\pi}{4})$
\end{tabular}\\ \hline

\begin{tabular}{c}
$SU(q+2) \curvearrowright $\\$Sp(q+2)/Sp(2)\times Sp(q)$\\$(q=3)$
\end{tabular}
&
\begin{tabular}{c}\\
$\includegraphics[height=50mm]{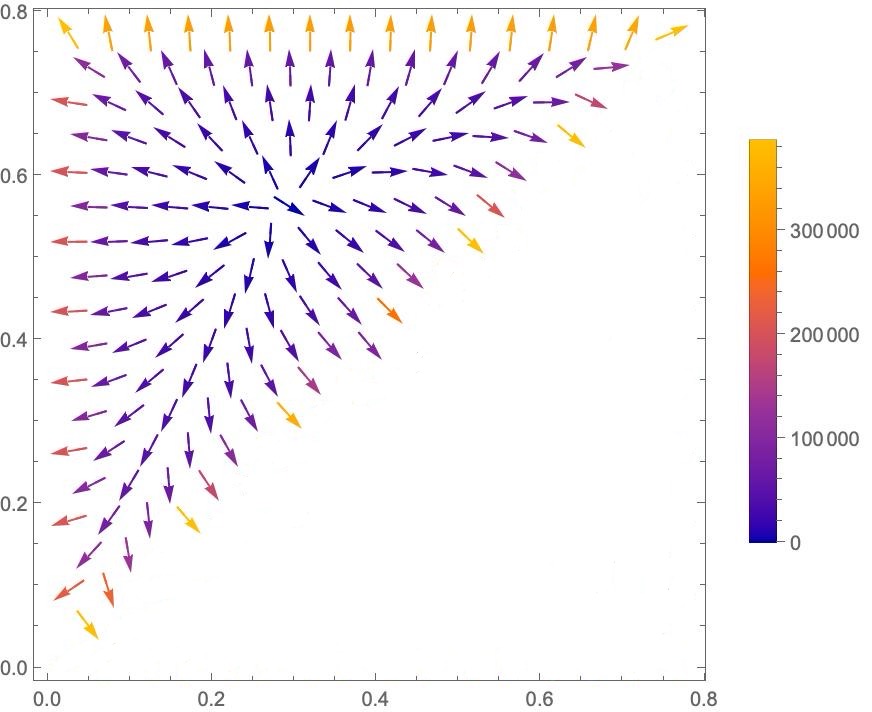}$\\
$(x_1^0,x_2^0) =$ (0.28431, 0.567081)\\
$(x_1^{\sigma},x_2^{\sigma}) = (0,0.53996),(0.42854,0.42854),(0.53996,\frac{\pi}{4})$
\end{tabular}\\ \hline


\hline
\begin{tabular}{c}
$SU(4)\curvearrowright $\\$Sp(4)/Sp(2)\times Sp(2)$
\end{tabular}
&
\begin{tabular}{c}\\
$\includegraphics[height=50mm]{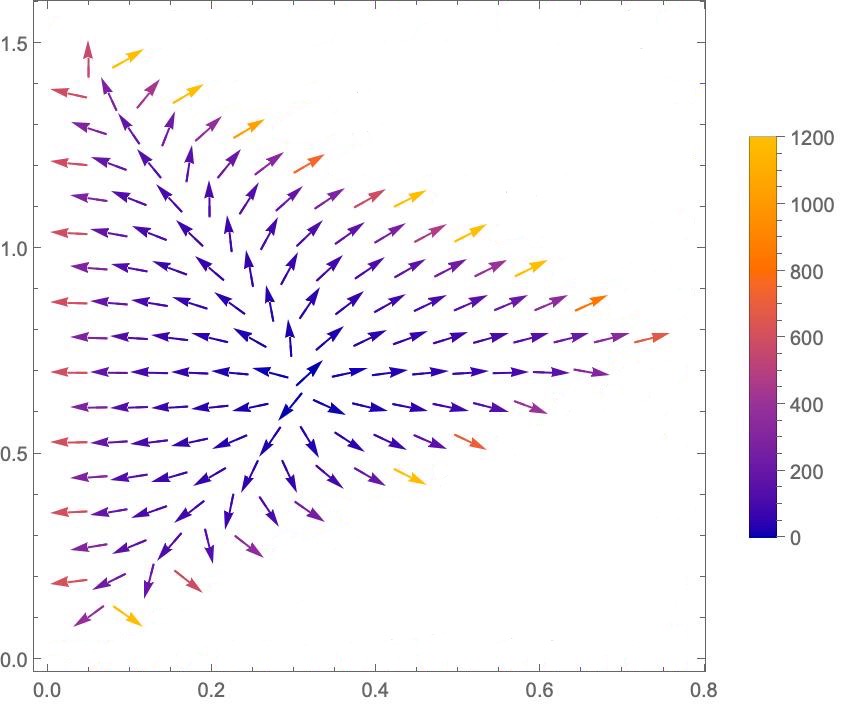}$\\
$(x_1^0,x_2^0) =$ (0.307799, 0.664173)\\
$(x_1^{\sigma},x_2^{\sigma}) = (0,0.68472), (0,0), (0,\frac{\pi}{2})$
\end{tabular}\\ \hline

\begin{tabular}{c}
$U(4)\curvearrowright $\\$Sp(4)/Sp(2)\times Sp(2)$
\end{tabular}
&
\begin{tabular}{c}\\
$\includegraphics[height=50mm]{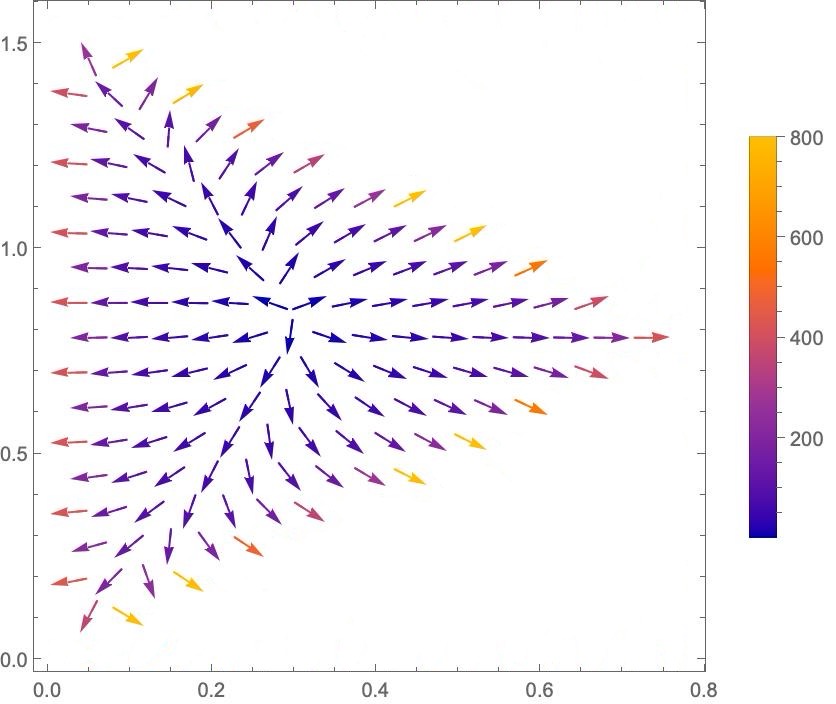}$\\
$(x_1^0,x_2^0) =$ (0.293247, 0.840516)\\
$(x_1^{\sigma},x_2^{\sigma}) = (0,0.88608), (0,0), (0,\frac{\pi}{2})$
\end{tabular}\\ \hline

\begin{tabular}{c}
$Sp(j+1)\times Sp(q-j+1)\curvearrowright $\\$(Sp(2)\times Sp(2))/Sp(2)$
\end{tabular}
&
\begin{tabular}{c}\\
$\includegraphics[height=50mm]{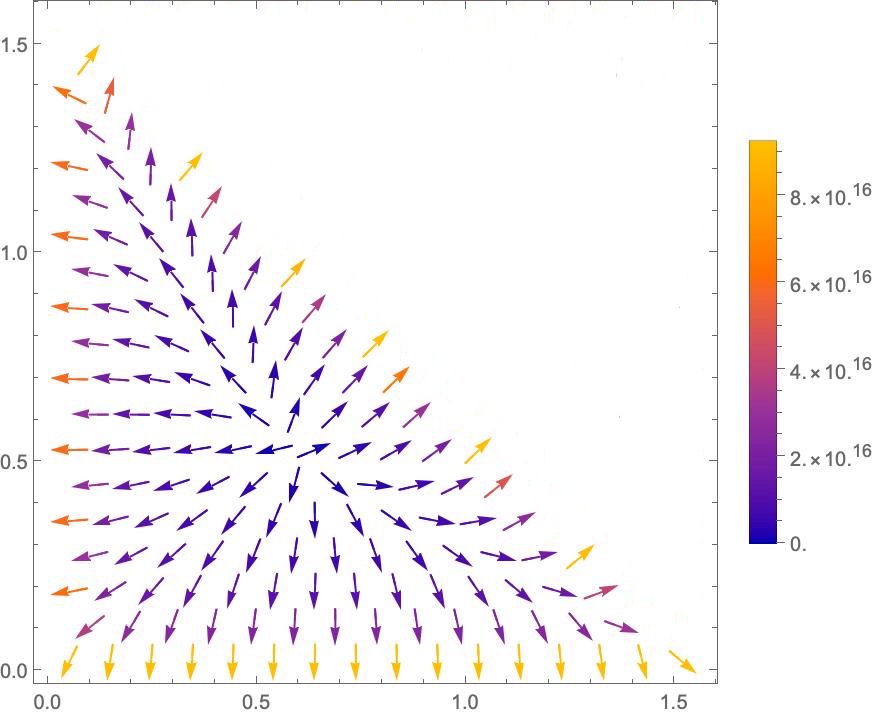}$\\
$(x_1^0,x_2^0) =$ (0.589609, 0.52777)\\
$(x_1^{\sigma},x_2^{\sigma}) = (0,0.42053), (0.67335,0), (0,\frac{\pi}{2})$
\end{tabular}\\ \hline


\hline
\begin{tabular}{c}
$Sp(2)\times Sp(2)\curvearrowright $\\$Sp(4)/Sp(2)\times Sp(2)$
\end{tabular}
&
\begin{tabular}{c}\\
$\includegraphics[height=50mm]{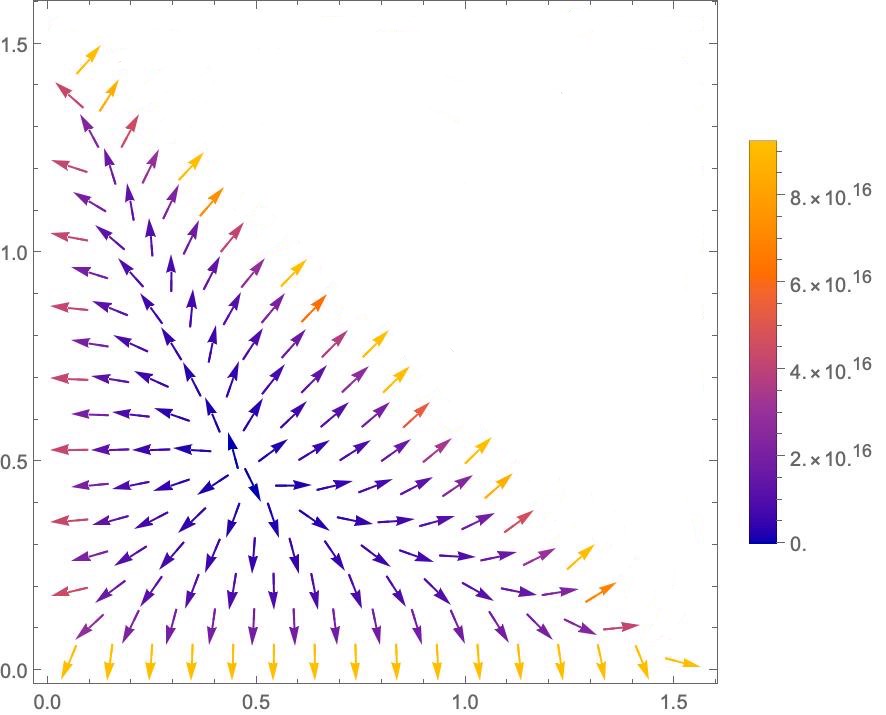}$\\
$(x_1^0,x_2^0) =$ (0.462616, 0.8711)\\
$(x_1^{\sigma},x_2^{\sigma}) = (0,0.57964), (0.54947,0), (0.930274,0.63355)$
\end{tabular}\\ \hline

\begin{tabular}{c}
$SU(2)^2\cdot SO(2)^2\curvearrowright $\\$(Sp(2)\times Sp(2))/Sp(2)$
\end{tabular}
&
\begin{tabular}{c}\\
$\includegraphics[height=50mm]{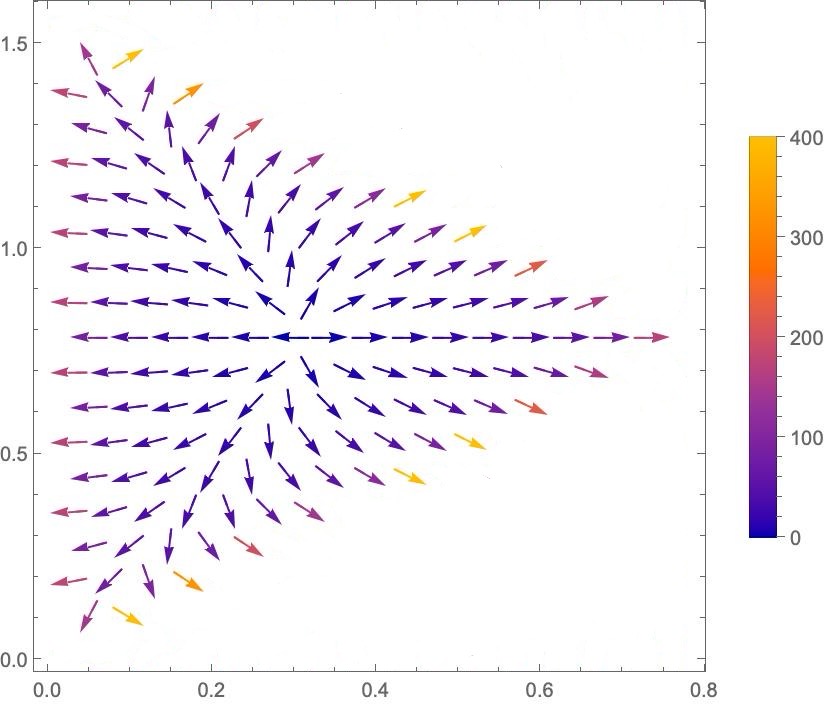}$\\
$(x_1^0,x_2^0) =$ (0.30774, 0.785398)\\
$(x_1^{\sigma},x_2^{\sigma}) = (0,\frac{\pi}{4}), (0,0), (0,\frac{\pi}{2})$
\end{tabular}\\ \hline

\begin{tabular}{c}
$\rho_{8}(Sp(2))\curvearrowright $\\$(Sp(2)\times Sp(2))/Sp(2)$
\end{tabular}
&
\begin{tabular}{c}\\
$\includegraphics[height=50mm]{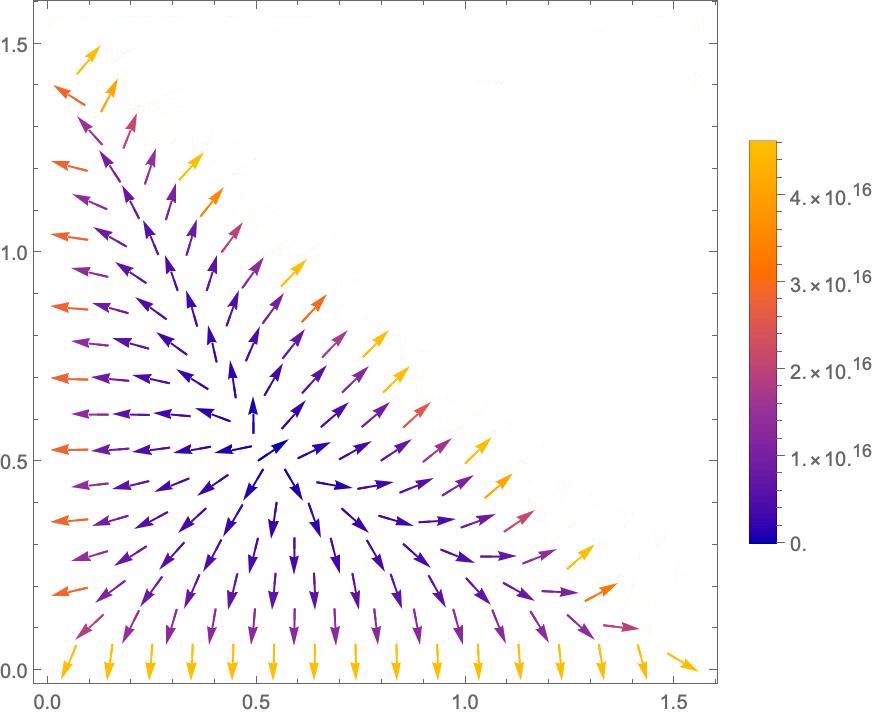}$\\
$(x_1^0,x_2^0) = (\frac{\pi}{6},\frac{\pi}{6})$
\end{tabular}\\ \hline


\hline
\begin{tabular}{c}
$\rho_{9}(Sp(2))\curvearrowright $\\$(Sp(2)\times Sp(2))/Sp(2)$
\end{tabular}
&
\begin{tabular}{c}\\
$\includegraphics[height=50mm]{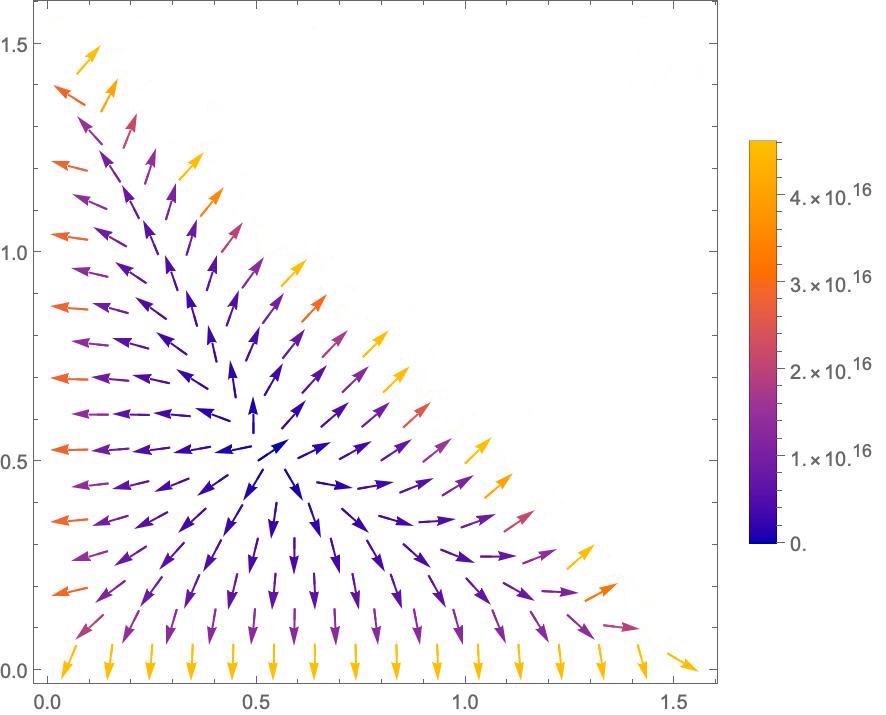}$\\
$(x_1^0,x_2^0) = (\frac{\pi}{6},\frac{\pi}{6})$
\end{tabular}\\ \hline

\begin{tabular}{c}
$Sp(4)\curvearrowright $\\$E_6/Spin(10)\cdot U(1)$
\end{tabular}
&
\begin{tabular}{c}\\
$\includegraphics[height=50mm]{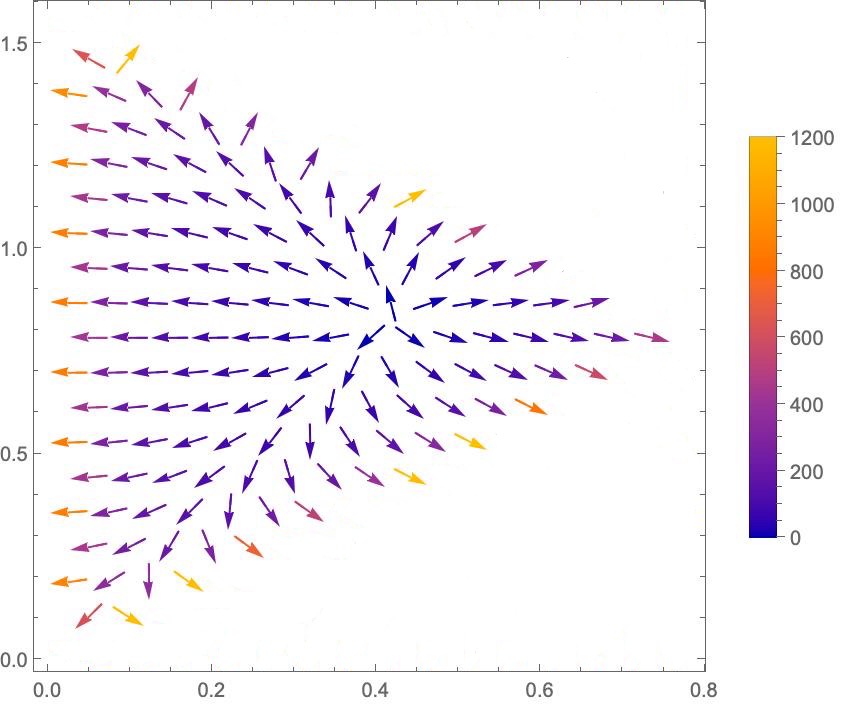}$\\
$(x_1^0,x_2^0) =$ (0.418547, 0.831616)\\
$(x_1^{\sigma},x_2^{\sigma}) = (0,\frac{\pi}{4}), (0,0), (0,\frac{\pi}{2})$
\end{tabular}\\ \hline

\begin{tabular}{c}
$SU(6)\cdot SU(2) \curvearrowright $\\$E_6/Spin(10)\cdot U(1)$
\end{tabular}
&
\begin{tabular}{c}\\
$\includegraphics[height=50mm]{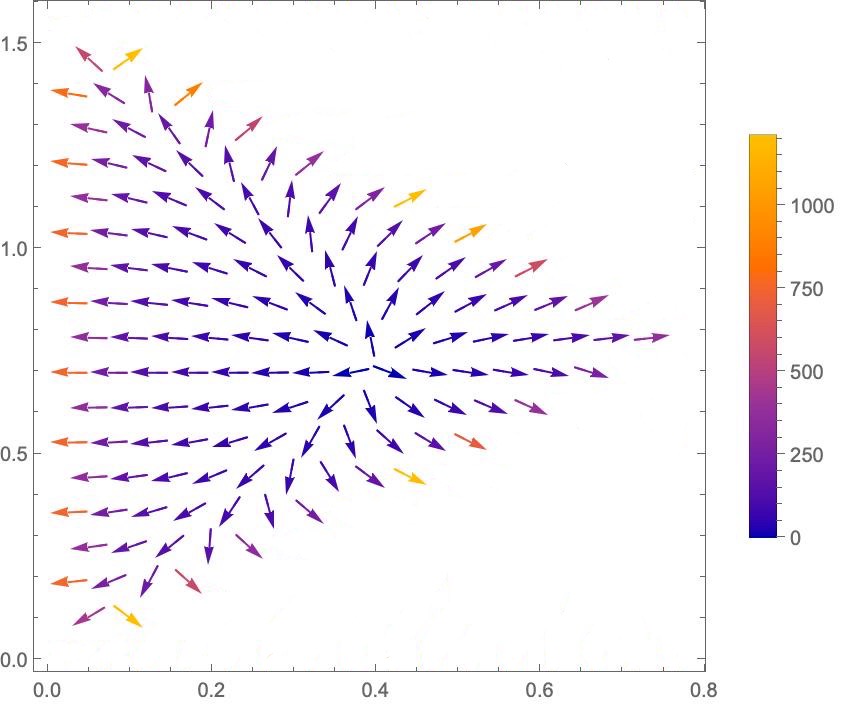}$\\
$(x_1^0,x_2^0) =$ (0.397099, 0.718451)\\
$(x_1^{\sigma},x_2^{\sigma}) = (0,\frac{\pi}{4}), (0,0), (0,\frac{\pi}{2})$
\end{tabular}\\ \hline


\hline
\begin{tabular}{c}
$\rho_{10}(SU(6)\cdot SU(2))\curvearrowright $\\$E_6/Spin(10)\cdot U(1)$
\end{tabular}
&
\begin{tabular}{c}\\
$\includegraphics[height=50mm]{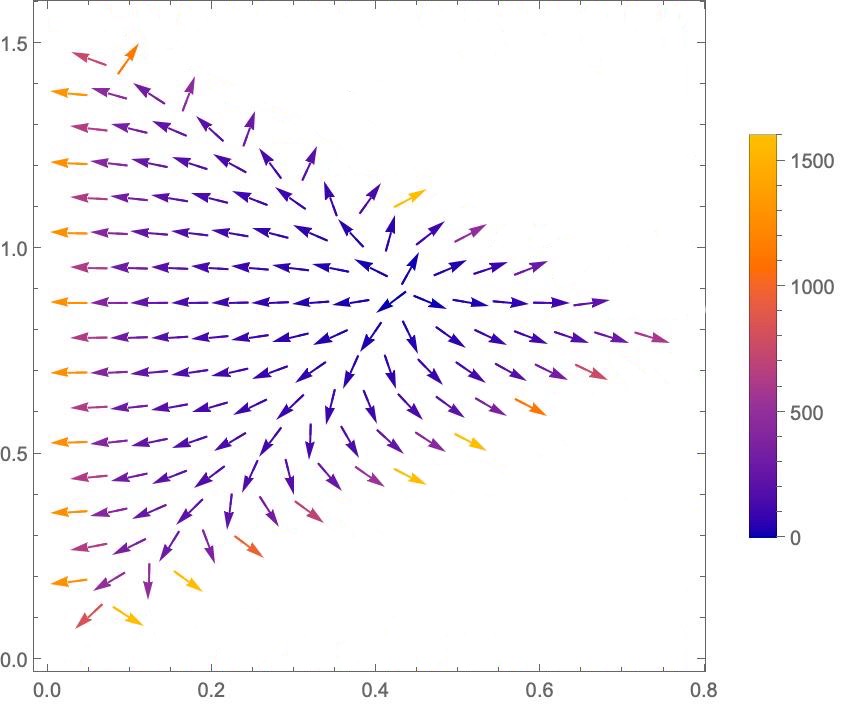}$\\
$(x_1^0,x_2^0) =$ (0.434443, 0.902558)
\end{tabular}\\ \hline

\begin{tabular}{c}
$\rho_{11}(Spin(10)\cdot SU(1))\curvearrowright $\\$E_6/Spin(10)\cdot U(1)$
\end{tabular}
&
\begin{tabular}{c}\\
$\includegraphics[height=50mm]{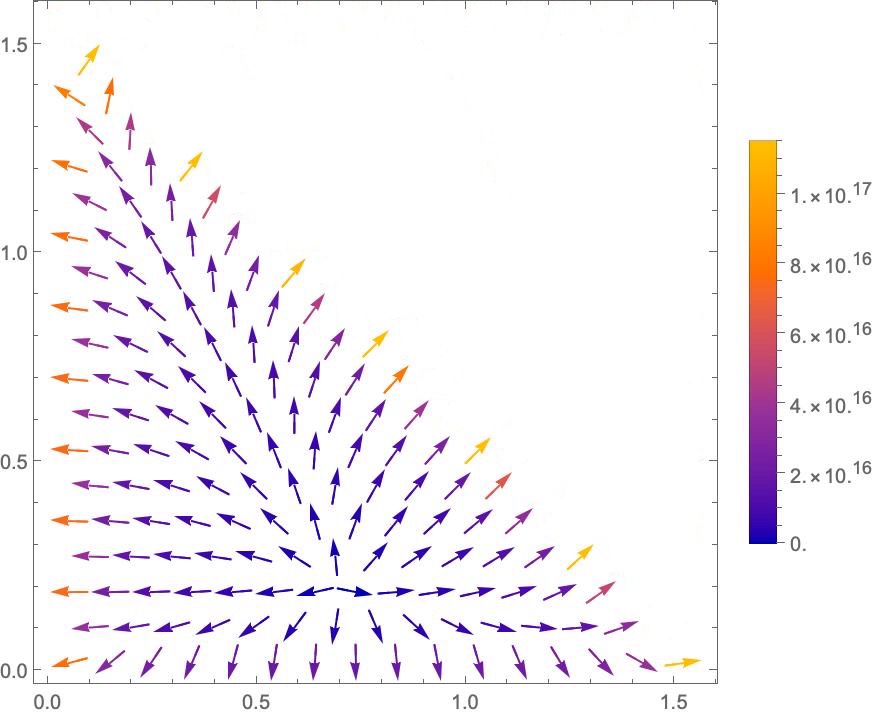}$\\
$(x_1^0,x_2^0) =$ (0.706305, 0.197372)
\end{tabular}\\ \hline

\begin{tabular}{c}
$\rho_{12}(Spin(10)\cdot SU(1))\curvearrowright $\\$E_6/Spin(10)\cdot U(1)$
\end{tabular}
&
\begin{tabular}{c}\\
$\includegraphics[height=50mm]{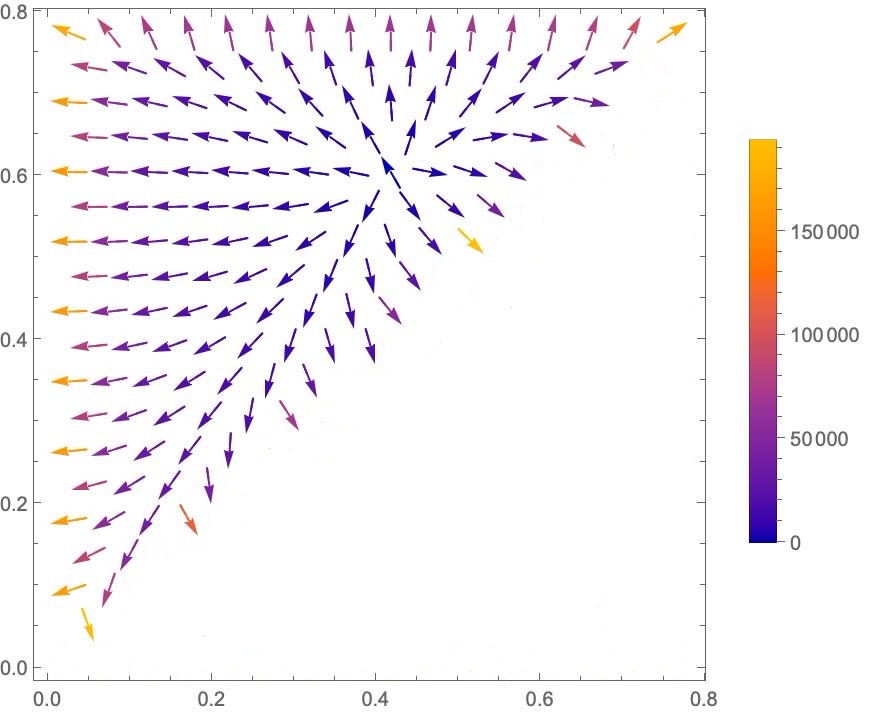}$\\
$(x_1^0,x_2^0) =$ (0.418704, 0.602223)
\end{tabular}\\ \hline


\hline
\begin{tabular}{c}
$Sp(4)\curvearrowright E_6/F_4$
\end{tabular}
&
\begin{tabular}{c}\\
$\includegraphics[height=50mm]{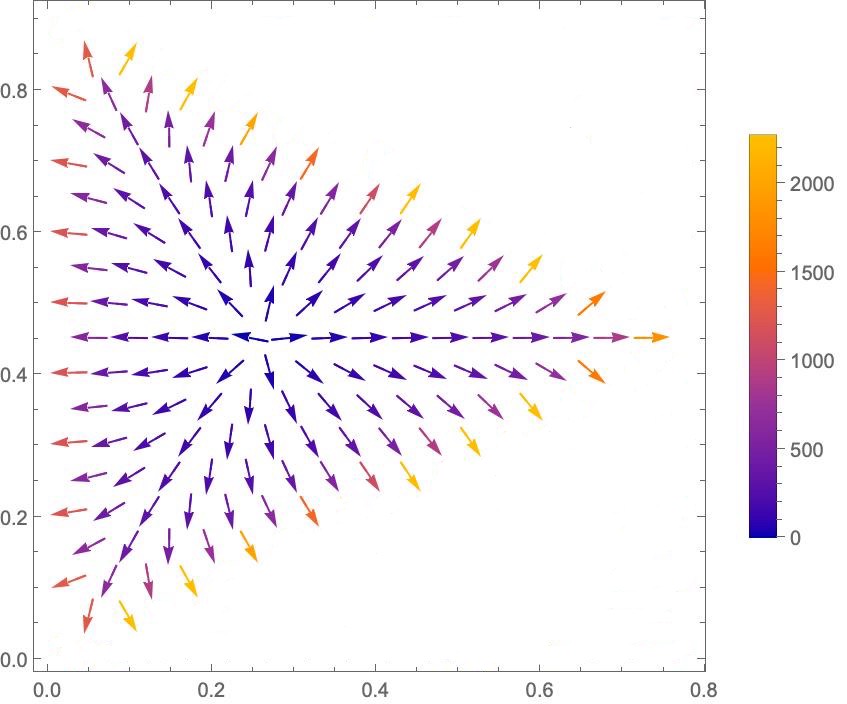}$\\
\end{tabular}\\ \hline

\begin{tabular}{c}
$\rho_{13}(F_4)\curvearrowright E_6/F_4$
\end{tabular}
&
\begin{tabular}{c}\\
$\includegraphics[height=50mm]{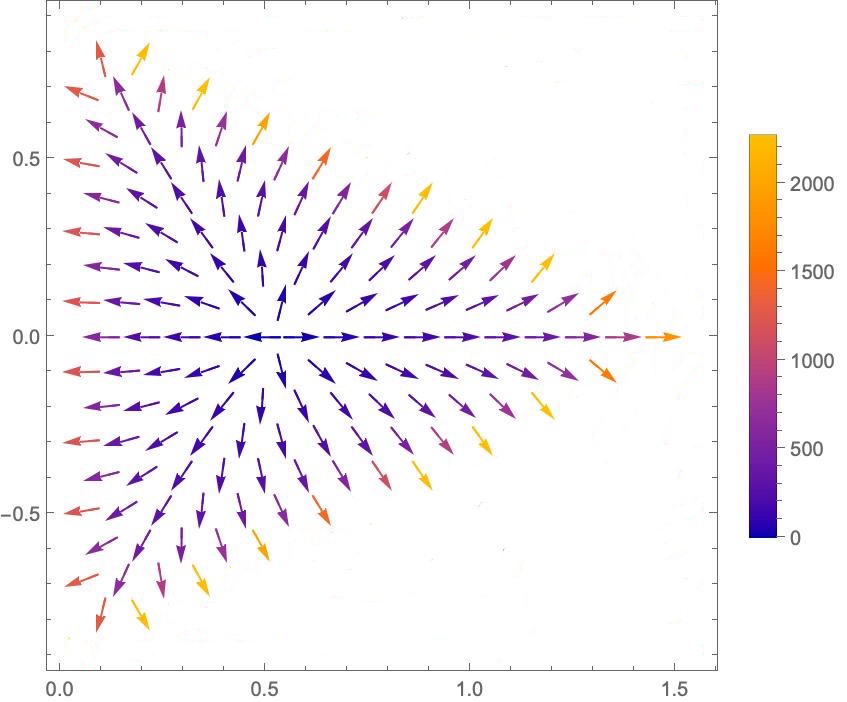}$\\
\end{tabular}\\ \hline

\begin{tabular}{c}
$\rho_{14}(SO(4))\curvearrowright  G_2/SO(4)$
\end{tabular}
&
\begin{tabular}{c}\\
$\includegraphics[height=50mm]{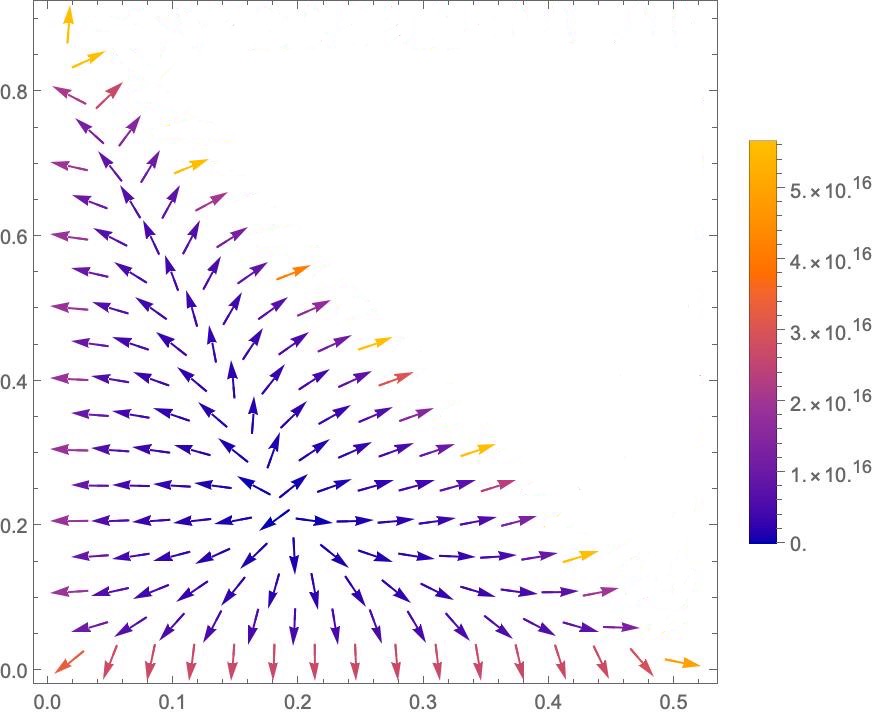}$\\
\end{tabular}\\ \hline


\hline
\begin{tabular}{c}
$\rho_{15}(SO(4))\curvearrowright  G_2/SO(4)$
\end{tabular}
&
\begin{tabular}{c}\\
$\includegraphics[height=50mm]{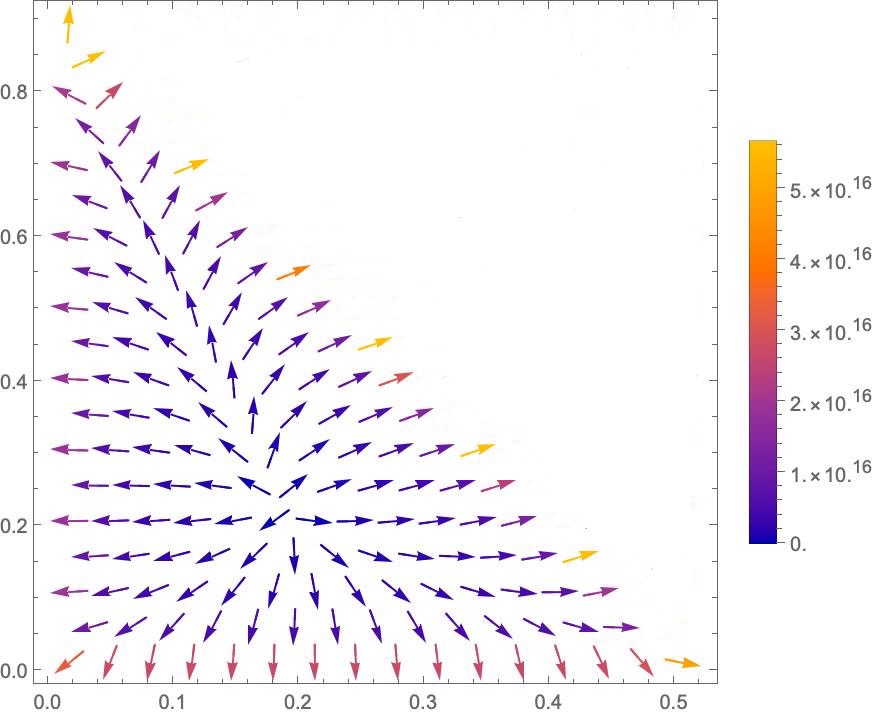}$\\
\end{tabular}\\ \hline

\begin{tabular}{c}
$\rho_{16}(G_2)\curvearrowright  (G_2\times G_2)/G_2$
\end{tabular}
&
\begin{tabular}{c}\\
$\includegraphics[height=50mm]{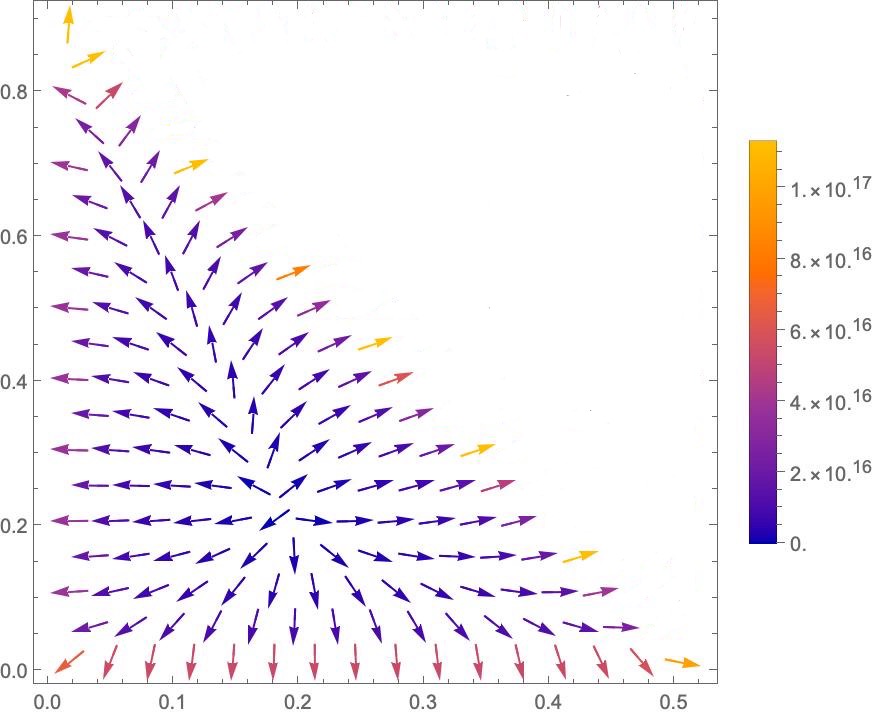}$\\
\end{tabular}\\ \hline

\begin{tabular}{c}
$SU(2)^4\curvearrowright  (G_2\times G_2)/G_2$
\end{tabular}
&
\begin{tabular}{c}\\
$\includegraphics[height=50mm]{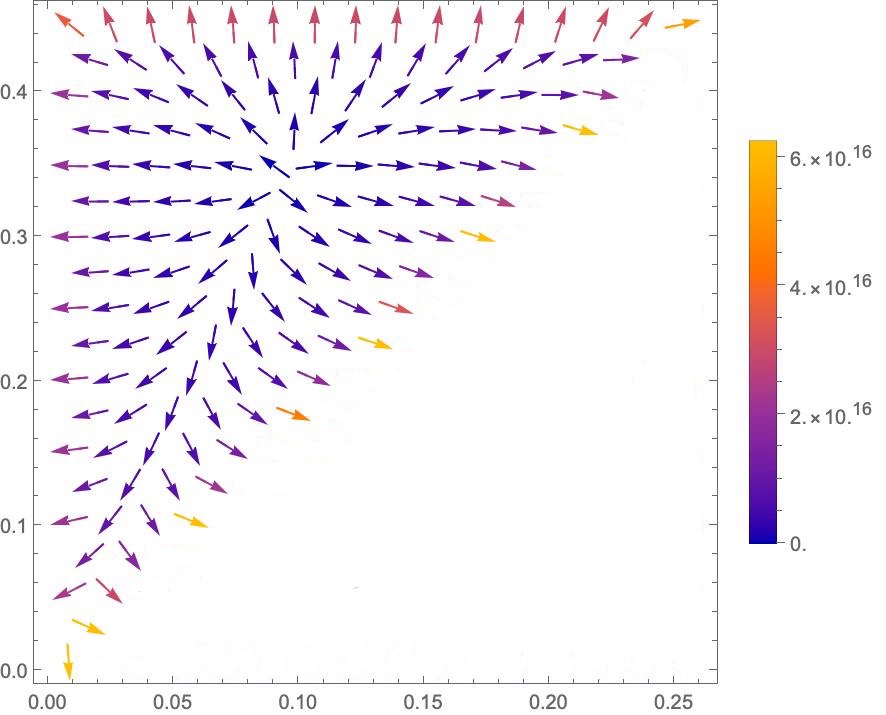}$\\
\end{tabular}\\ \hline


\end{longtable}

\vspace{0.5truecm}

{\small 
\rightline{Department of Mathematics, Faculty of Science}
\rightline{Tokyo University of Science, 1-3 Kagurazaka}
\rightline{Shinjuku-ku, Tokyo 162-8601 Japan}
\rightline{(koike@rs.tus.ac.jp)}
}

\vspace{0.5truecm}

{\small 
\rightline{Department of Mathematics, Faculty of Science}
\rightline{Tokyo University of Science, 1-3 Kagurazaka}
\rightline{Shinjuku-ku, Tokyo 162-8601 Japan}
\rightline{(1123514@ed.tus.ac.jp)}
}


\vspace{0.5truecm}

\begin{thebibliography}{999}

\bibitem[HPTT]{HPTT}  E. Heintze, R. S. Palais, C. L. Terng, and G. Thorbergsson, Hyperpolar actions on symmetric spaces, 
Geometry, topology and physics for Raoul Bott, (ed. S.T.Yau), Conf. Proc. Lecture Notes  Geom. Topology 4, Internat. Press, 
Cambridge, MA, (1995), 214-245.
\bibitem[Koi]{Koi}  
N. Koike, Collapse of the mean curvature flow for equifocal submanifolds, Asian J. Math {\bf 15} (2011), 101--128.
\bibitem[Kol]{Kol}  
A. Kollross, A Classification of hyperpolar and cohomogeneity one actions, Trans. Amer. Math. Soc., {\bf 354} (2001), 571--612.
\bibitem[LT]{LT}  
X. Liu and C. L. Terng, The mean curvature flow for isoparametric submanifolds, Duke Math. J. {\bf 147} (2009), 157--179.
\bibitem[PT]{PT} R. S. Palais and C. L. Terng, Critical point theory and submanifold geometry, Lecture Notes 
in Math. {\bf 1353}, Springer, Berlin, 1988.
\bibitem[TT]{TT}  
C.L. Terng, and G. Thorbergsson, Submanifold geometry in symmetric spaces, J. Differential Geom. {\bf 42} (1995), 665--718.

\end{thebibliography}
\end{document}